\documentclass[a4paper,reqno,12pt]{amsart}
\usepackage{amsmath,amsrefs,color}

\numberwithin{equation}{section}

\DeclareMathOperator{\meas}{meas}
\DeclareMathOperator{\loc}{loc}
\newcommand{\OO}{\text{O}}
\newcommand{\oo}{\text{o}}
\newcommand{\R}{\mathbb{R}}
\newcommand{\N}{\mathbb{N}}

\renewcommand{\>}{\right>}
\renewcommand{\[}{\left[}
\renewcommand{\]}{\right]}
\renewcommand{\(}{\left(}
\renewcommand{\)}{\right)}

\newtheorem{theorem}{Theorem}[section]
\newtheorem{lemma}[theorem]{Lemma}
\newtheorem{step}[theorem]{Step}

\begin{document}

\title[Finite energy solutions to elliptic systems in $\R^n$]{Decay estimates and symmetry of finite energy solutions to elliptic systems in $\R^n$}

\author{J\'er\^ome V\'etois}

\address{J\'er\^ome V\'etois, McGill University Department of Mathematics and Statistics, 805 Sherbrooke Street West, Montreal, Quebec H3A 0B9, Canada.}
\email{jerome.vetois@mcgill.ca}

\date{December 29, 2016}

\begin{abstract}
We study a notion of finite energy solutions to elliptic systems with power nonlinearities in $\R^n$. We establish sharp pointwise decay estimates for positive and sign-changing solutions. By using these estimates, we obtain symmetry results when the solutions are positive.
\end{abstract}

\maketitle

\section{Introduction and main results}\label{Sec1}

In this paper, we are interested in the nonlinear Schr\"odinger system
\begin{equation}\label{Eq1}
\left\{\begin{aligned}&-\Delta u=\left|v\right|^{p}\left|u\right|^{r-1}u&&\text{in }\R^n\\
&-\Delta v=\left|u\right|^{q}\left|v\right|^{s-1}v&&\text{in }\R^n
\end{aligned}\right.
\end{equation}
where 
\begin{equation}\label{Eq2}
n\ge3,\quad p,q\ge1,\quad r,s\ge0,\quad\text{and}\quad p-s\ge q-r>-1.
\end{equation}
Systems of this type are closely related to the modeling of Bose--Einstein condensates (see for instance Burke, Bohn, Esry and Greene~\cite{BurBohnEsryGre} and Frantzeskakis~\cite{Fra}) and also arise in nonlinear optics (see Kanna and Lakshmanan~\cite{KanLak} and Kivshar and Luther-Davies~\cite{KivLut}).

\smallskip
We define 
$$a:=\frac{n\(pq-\(r-1\)\(s-1\)\)}{2\(p-s+1\)}\quad\text{and}\quad b:=\frac{n\(pq-\(r-1\)\(s-1\)\)}{2\(q-r+1\)}$$
so that
\begin{equation}\label{Eq3}
\frac{p}{b}+\frac{r}{a}=\frac{2}{n}+\frac{1}{a}\quad\text{and}\quad\frac{q}{a}+\frac{s}{b}=\frac{2}{n}+\frac{1}{b}\,.
\end{equation}
It follows from \eqref{Eq3} that \eqref{Eq1} is invariant under the change of scale
\begin{equation}\label{Eq4}
\(\mu^{n/a}u\(\mu x\),\mu^{n/b}v\(\mu x\)\)
\end{equation}
for all $\mu>0$. Moreover, the $L^a$--norm of $u$ and $L^b$--norm of $v$ are also invariant under \eqref{Eq4}. We assume that
\begin{equation}\label{Eq5}
a,b>\frac{n}{n-2}\,.
\end{equation}
It follows from \eqref{Eq3} and \eqref{Eq5} that for every $\(u,v\)\in L^a\(\R^n\)\times L^b\(\R^n\)$, $\left|v\right|^{p}\left|u\right|^{r-1}u\in L^1_{\loc}\(\R^n\)$ and $\left|u\right|^{q}\left|v\right|^{s-1}v\in L^1_{\loc}\(\R^n\)$. We then say that $\(u,v\)\in L^a\(\R^n\)\times L^b\(\R^n\)$ is a {\it solution} (or {\it finite energy solution}) of \eqref{Eq1} if both equations in \eqref{Eq1} are satisfied in the distributional sense. 

\smallskip
By using Sobolev's inequalities, it is easy to check that in case $r=s=0$, if $a$ and $b$ satisfy the critical condition
\begin{equation}\label{Eq6}
\frac{1}{a}+\frac{1}{b}=\frac{n-2}{n}\,,
\end{equation}
then our definition of finite energy solutions generalizes the more usual definition of finite energy solutions as critical points in Sobolev spaces of an energy functional (see Lions~\cite{Lio} for the variational framework in this case, see also Hulshof and van der Vorst~\cite{HulVan}). Remark that this definition cannot be directly extended since the system \eqref{Eq1} is not variational in the general case where $r$ and $s$ can be nonzero. In this case, our definition provides a natural generalization of such solutions which preserves the scale invariance \eqref{Eq4}.

\smallskip\addtolength{\textheight}{-14pt}
It follows from a result of Chen, Li, and Ou~\cite{ChenLiOu} that in case $r=s=0$, if \eqref{Eq6} holds true, then positive finite energy solutions are radially symmetric and decreasing about some point $x_0\in\R^n$. Moreover, positive finite energy solutions have been proven to exist in this case if and only if \eqref{Eq6} holds true (see Lions~\cite{Lio} for the existence and Chen and Li~\cite{ChenLi2} and Villavert~\cite{Vil} for the nonexistence). We also mention that positive solutions with infinite energy of \eqref{Eq1} have been proven to exist in case $r=s=0$ when the equality in \eqref{Eq6} is replaced by $<$ (see Mitidieri~\cite{Mit} and Serrin and Zou~\cite{SerZou3}). In case of the converse inequality $>$, nonexistence of positive solutions (with finite or infinite energy) of \eqref{Eq1} with $r=s=0$ has been obtained in the radial case (see \nobreak Mitidieri~\cite{Mit} and Serrin and Zou~\cite{SerZou1}), in case $n\le4$ (see Mitidieri~\cite{Mit}, Serrin and Zou~\cite{SerZou2}, and Souto~\cite{Souto} in case $n\in\left\{1,2\right\}$, Pol\'a\v{c}ik, Quittner, and \nobreak  Souplet~\cite{PolQuiSou} and Serrin and Zou~\cite{SerZou2} in case $n=3$, and Souplet~\cite{Soup} in case $n=4$), and under various additional assumptions in higher dimensions (see for instance Busca and Man\'asevich~\cite{BusMan}, Chen and Li~\cite{ChenLi2}, de Figueiredo and Felmer~\cite{deFFel}, Lei and Li~\cite{LeiLi2}, Quittner, and Souplet~\cite{PolQuiSou}, Mitidieri~\cite{Mit}, Serrin and Zou~\cites{SerZou1,SerZou2}, Souplet~\cite{Soup}, and Villavert~\cite{Vil}). We also refer to Lei and Li~\cite{LeiLi1} and Villavert~\cite{Vil} for a qualitative analysis of solutions of \eqref{Eq1} in case $r=s=0$.

\smallskip
Another case which has received a lot of attention is the case $q+s=p+r$. In this case, Quittner and Souplet~\cite{QuiSou} obtained that every positive solution $\(u,v\)$ of \eqref{Eq1} is such that either $u\equiv0$, $v\equiv0$, or $u\equiv v$ in $\R^n$ provided $p\ge r$ and $r,s\le n/\(n-2\)$ (see also Montaru, Sirakov, and Souplet~\cite{MonSirSou} for extensions of this result and Guo and Liu~\cite{GuoLiu} and Li and Ma~\cite{LiMa} for previous results in case $p+r=\(n+2\)/\(n-2\)$). In case $u\equiv v$, we are then left with the equation 
$$-\Delta u=u^{p+r-1}\qquad\text{in }\R^n.$$
In particular, we can then apply Gidas and Spruck's nonexistence result~\cite{GidSpr} in case $p+r<\(n+2\)/\(n-2\)$ and Caffarelli, Gidas and Spruck's classification result~\cite{CafGidSpr} in case $p+r=\(n+2\)/\(n-2\)$ (see also Chen and Li~\cite{ChenLi1} for an easier proof and Gidas, Ni, Nirenberg~\cite{GidNiNir} and Obata~\cite{Oba} for previous results). 

\smallskip
In the general case where $r$ and $s$ can be nonzero and $p+r$ can be different to $q+s$, nonexistence results of positive solutions of \eqref{Eq1} have been obtained under various subcritical-type conditions (see for instance Bidaut-V\'eron and Giacomini~\cite{BidVerGia}, Bidaut-V\'eron and Pohozaev \cite{BidVerPoh}, Chen and Lu~\cite{ChenLu}, Cl{\'e}ment, Fleckinger, Mitidieri, and de Th{\'e}lin~\cite{CleFleMitdeT}, Mitidieri~\cite{Mit}, Quittner~\cite{Qui}, Reichel and Zou~\cite{ReiZou}, and Zheng~\cite{Zhe}). Existence of radially symmetric solutions has also been established under various supercritical-type conditions (see for instance Bidaut-V\'eron and Giacomini~\cite{BidVerGia} and Li and Villavert~\cites{LiVil1,LiVil2}). However the optimal conditions for existence of positive solutions of \eqref{Eq1} are not yet known even for radially symmetric solutions.

\smallskip
Note that in contrast with the case of a single equation (see Caffarelli, Gidas and Spruck~\cite{CafGidSpr} and Chen and Li~\cite{ChenLi1}), the system \eqref{Eq1} is not invariant under the Kelvin transform in the critical case \eqref{Eq6} and the moving plane method does not seem to be directly applicable to the transformed system when $p+r>\(n+2\)/\(n-2\)$. 

\smallskip
In Theorem~\ref{Th1} below, we establish sharp pointwise decay estimates for finite energy solutions (positive or not) of \eqref{Eq1}. We then use these estimates to obtain symmetry results for positive solutions of \eqref{Eq1} in Theorems~\ref{Th2}--\ref{Th4}. These latter results extend previous works of Liu and Ma~\cites{LiuMa1,LiuMa2} who obtained symmetry results under decay assumptions on the solutions.

\smallskip
With regard to the signs of $u$ and $v$, we assume that 
\begin{equation}\label{Eq8}
u\ge0\text{ in }\R^n\text{ if }r=0\quad\text{and}\quad v\ge0\text{ in }\R^n\text{ if }s=0
\end{equation}
so that the nonlinearities in \eqref{Eq1} are continuous functions of $u$ and $v$. Moreover, we assume that
\begin{equation}\label{Eq9}
u\ge0\text{ in }\R^n\text{ if }r<1\text{ and }s<1.
\end{equation}
The assumption \eqref{Eq9} allows us to use a comparison result between $u$ and $v$ in the sublinear case $r,s<1$ (see Step~\ref{Th1St2}).

\smallskip
We obtain the following decay estimates:

\begin{theorem}\label{Th1}
Assume that \eqref{Eq2} and \eqref{Eq5} hold true. Let $\(u,v\)\in L^a\(\R^n\)\times L^b\(\R^n\)$ be a solution of \eqref{Eq1} such that \eqref{Eq8} and \eqref{Eq9} hold true. Then $u,v\in C^2\(\R^n\)$ and there exists a constant $C_0$ such that
\begin{equation}\label{Th1Eq1}
\left\{\begin{aligned}&\left|u\(x\)\right|\le C_0(1+\left|x\right|^{n-2})^{-1}\\&\left|v\(x\)\right|\le C_0\(1+h\(x\)\)^{-1}\end{aligned}\right.\qquad\forall x\in\R^n
\end{equation}
where
\begin{equation}\label{Th1Eq2}
h\(x\):=\left\{\begin{aligned}&\left|x\right|^{n-2}&&\text{if }q+s>n/\(n-2\)\\&\left|x\right|^{n-2}\ln\(1+\left|x\right|\)^{-1/\(1-s\)}&&\text{if }q+s=n/\(n-2\)\\&\left|x\right|^{\(\(n-2\)q-2\)/\(1-s\)}&&\text{if }q+s<n/\(n-2\).\end{aligned}\right.
\end{equation}
Moreover, if $u,v\ge0$ in $\R^n$, then either $u\equiv v\equiv0$ in $\R^n$ or there exists a constant $C_1>0$ such that
\begin{equation}\label{Th1Eq3}
\left\{\begin{aligned}&u\(x\)\ge C_1(1+\left|x\right|^{n-2})^{-1}\\&v\(x\)\ge C_1\(1+h\(x\)\)^{-1}\end{aligned}\right.\qquad\forall x\in\R^n.
\end{equation}
\end{theorem}

Note that by using \eqref{Eq2}--\eqref{Eq5}, we obtain that if $q+s\le n/\(n-2\)$, then $s<1$ and $q>2/\(n-2\)$ hence the exponents in \eqref{Th1Eq2} are well defined and positive.

\smallskip
In case of the equation $-\Delta u=\left|u\right|^{4/\(n-2\)}u$ (i.e. $p+r=q+s=\(n+2\)/\(n-2\)$ and $u\equiv v$ in $\R^n$), the estimate \eqref{Th1Eq1} was obtained by Jannelli and Solimini~\cite{JanSol} (see also V\'etois~\cite{Vet} and Xiang~\cite{Xia} for different proofs of this result and generalizations to $p$--Laplace equations).

\smallskip
When $u,v\ge0$ in $\R^n$ and $q+s\ge n/\(n-2\)$, we obtain the following symmetry result:

\begin{theorem}\label{Th2}
Assume that \eqref{Eq2} and \eqref{Eq5} hold true. Assume moreover that $q+s\ge n/\(n-2\)$. Then for any solution $\(u,v\)\in L^a\(\R^n\)\times L^b\(\R^n\)$ of \eqref{Eq1} such that $u,v\ge0$ in $\R^n$, the functions $u$ and $v$ are radially symmetric and decreasing about some point $x_0\in\R^n$.
\end{theorem}

We point out that in case $q+s>n/\(n-2\)$, the result of Theorem~\ref{Th2} can be obtained directly by combining Theorem~\ref{Th1} and results of Liu and Ma~\cites{LiuMa1,LiuMa2}.

\smallskip
In case $q+s<n/\(n-2\)$, we obtain the following result:

\begin{theorem}\label{Th3}
Assume that \eqref{Eq2} and \eqref{Eq5} hold true. Assume moreover that $q+s<n/\(n-2\)$. Let $\(u,v\)\in L^a\(\R^n\)\times L^b\(\R^n\)$ be a solution of \eqref{Eq1} such that $u,v\ge0$ in $\R^n$ and
\begin{equation}\label{Th3Eq}
\ell_u^q\ell_v^{s-1}<C_{n,q,s}
\end{equation}
where 
$$\ell_u:=\limsup_{\left|x\right|\to\infty}\(\left|x\right|^{n-2}u\(x\)\)\text{ and }\ell_v:=\liminf_{\left|x\right|\to\infty}\big(\left|x\right|^{\(q\(n-2\)-2\)/\(1-s\)}v\(x\)\big)$$
and 
$$C_{n,q,s}:=\left\{\begin{aligned}&\frac{\(n-2\)^2}{4s}&&\text{if }2q+s\ge\frac{n+2}{n-2}\\&\frac{\(\(n-2\)q-2\)\(n-\(n-2\)\(q+s\)\)}{s\(1-s\)^2}&&\text{if }2q+s<\frac{n+2}{n-2}\end{aligned}\right.$$
in case $s\in\(0,1\)$ and $C_{n,q,s}=\infty$ in case $s=0$ (hence \eqref{Th3Eq} is always true in case $s=0$). Then the functions $u$ and $v$ are radially symmetric and decreasing about some point $x_0\in\R^n$.
\end{theorem}

Remark that it follows from Theorem~\ref{Th1} that we always have $\ell_u<\infty$ and $\ell_v>0$ when $q+s<n/\(n-2\)$. In case $u$ and $v$ are asymptotically equivalent to power functions, namely such that
\begin{equation}\label{Eq10}
\lim_{\left|x\right|\to\infty}\(\left|x\right|^{n-2}u\(x\)\)=\ell_u
\end{equation}
and
\begin{equation}\label{Eq11}
\lim_{\left|x\right|\to\infty}\big(\left|x\right|^{\(q\(n-2\)-2\)/\(1-s\)}v\(x\)\big)=\ell_v,
\end{equation}
we obtain the following result:

\begin{theorem}\label{Th4}
Assume that \eqref{Eq2} and \eqref{Eq5} hold true. Assume moreover that $q+s<n/\(n-2\)$. Let $\(u,v\)\in L^a\(\R^n\)\times L^b\(\R^n\)$ be a solution of \eqref{Eq1} such that $u,v\ge0$ in $\R^n$ and \eqref{Eq10} and \eqref{Eq11} hold true. Then 
$$\ell_u^q\ell_v^{s-1}=\frac{\(\(n-2\)q-2\)\(n-\(n-2\)\(q+s\)\)}{\(1-s\)^2}$$
which is smaller than $C_{n,q,s}$ hence it follows from Theorem~\ref{Th3} that the functions $u$ and $v$ are radially symmetric and decreasing about some point $x_0\in\R^n$.
\end{theorem}

It is well-known that solutions of Schr\"odinger equations and systems in $\R^n$ like \eqref{Eq1} play a fundamental role in the blow-up analysis of solutions of more general equations on arbitrary domains or manifolds. We refer to Struwe~\cite{Str} as an historic reference on this topic in case of the equation $-\Delta u=\left|u\right|^{4/\(n-2\)}u$. In case of systems, some references on this topic are for instance Druet and Hebey~\cite{DruHeb}, Hebey~\cites{Heb1,Heb2}, and Thizy~\cite{Thi}. With regard to the study of finite energy solutions for systems in $\R^n$, it is also interesting to mention the recent works of Gladiali, Grossi, and Troestler~\cite{GlaGroTro} and Guo, Li, and Wei~\cite{GuoLiWei} where existence of nonradial finite energy solutions is obtained for a class of systems with sums of power nonlinearities.

\smallskip
The proof of Theorem~\ref{Th1} is in two parts. First in Section~\ref{Sec1}, we obtain regularity and integrability results in weak Lebesgue spaces. Then in Section~\ref{Sec2}, we use the results of Section~\ref{Sec1} to obtain sharp pointwise estimates. The proof of Theorem~\ref{Th1} also relies on two preliminary steps which are on the one hand, a preliminary pointwise estimate that we obtain by using blow-up arguments and a doubling property (see Pol\'a\v{c}ik, Quittner, and Souplet~\cite{PolQuiSou}) and on the other, a generalization of comparison results obtained by Quittner and Souplet~\cite{QuiSou} in case $p+r=q+s$ and Souplet~\cite{Soup} in case $r=s=0$. We prove Theorems~\ref{Th2}--\ref{Th4} in Section~\ref{Sec3}. The proof of Theorems~\ref{Th2} and~\ref{Th3} is based on the moving plane method. As in the papers of Liu and Ma~\cites{LiuMa1,LiuMa2} (see also Chen and Li~\cites{ChenLi1,ChenLi2}), we apply the moving plane method to auxiliary functions with lower decay rates at infinity than the functions $u$ and $v$. We manage to extend the proof of Liu and Ma~\cites{LiuMa1,LiuMa2} in case $q+s\le n/\(n-2\)$ by allowing our two auxiliary functions to have different decay rates at infinity.

\smallskip\noindent
{\bf Acknowledgments.} The author was supported by a Discovery Grant awarded by the Natural Sciences and Engineering Research Council of Canada. The author is very grateful to Mostafa Fazly for fruitful discussions and several references on systems of type \eqref{Eq1}. The author is also very grateful to Pavol Quittner and Frederic Robert for helpful comments on earlier versions of this paper.

\section{Regularity and integrability results}\label{Sec2}

For any Lebesgue measurable set $\Omega\subset\R^n$ and any $\sigma\in\(0,\infty\)$ and $\tau\in\(0,\infty\]$, we let $L^{\sigma,\tau}\(\Omega\)$ be the Lorentz space defined as the set of all measurable functions $f:\Omega\to\R$ such that $\left\|f\right\|_{L^{\sigma,\tau}\(\Omega\)}<\infty$, where
$$\left\|f\right\|_{L^{\sigma,\tau}\(\Omega\)}:=\left\{\begin{aligned}&\sigma^{1/\tau}\(\int_0^\infty h^{\tau-1}\meas\(\left\{\left|f\right|>h\right\}\)^{\tau/\sigma}dh\)^{1/\tau}&&\text{if }\tau<\infty\\&\sup_{h>0}\big(h^\sigma\meas\(\left\{\left|f\right|>h\right\}\)\big)^{1/\sigma}&&\text{if }\tau=\infty\end{aligned}\right.$$
with $\meas\(\left\{\left|f\right|>h\right\}\)$ standing for the Lebesgue measure of the set $\left\{x\in\Omega:\,\left|f\(x\)\right|>h\right\}$. We refer to the book of Grafakos~\cite{Gra}*{Chapter~1} for a general presentation of Lorentz spaces.

\smallskip

In Lemmas~\ref{Lem1},~\ref{Lem2}, and~\ref{Lem3} below, we recall the generalizations of H\"older's and Young's inequalities for Lorentz spaces. These results are due to O'Neil~\cite{ONe}.

\begin{lemma}{\rm(H\"older's inequality)}\label{Lem1}
Let $\Omega$ be a Lebesgue measurable subset of $\R^n$ and $\sigma,\sigma_1,\sigma_2\in\(0,\infty\)$ and $\tau,\tau_1,\tau_2\in\(0,\infty\]$ be such that 
\begin{equation}\label{Lem1Eq}
\frac{1}{\sigma_1}+\frac{1}{\sigma_2}=\frac{1}{\sigma}\quad\text{and}\quad\frac{1}{\tau_1}+\frac{1}{\tau_2}\ge\frac{1}{\tau}
\end{equation}
with the convention that $1/\infty=0$. Then for any $f_1\in L^{\sigma_1,\tau_1}\(\Omega\)$ and $f_2\in L^{\sigma_2,\tau_2}\(\Omega\)$, we have $f_1f_2\in L^{\sigma,\tau}\(\Omega\)$ and
$$\left\|f_1f_2\right\|_{L^{\sigma,\tau}\(\Omega\)}\le C\left\|f_1\right\|_{L^{\sigma_1,\tau_1}\(\Omega\)}\left\|f_2\right\|_{L^{\sigma_2,\tau_2}\(\Omega\)}$$
for some constant $C$ independent of $\Omega$, $f_1$, and $f_2$.
\end{lemma}

\begin{lemma}{\rm(Young's inequality)}\label{Lem2}
Let $\sigma,\sigma_1,\sigma_2\in\(1,\infty\)$, and $\tau,\tau_1,\tau_2\in\(0,\infty\]$ be such that 
\begin{equation}\label{Lem2Eq}
\frac{1}{\sigma_1}+\frac{1}{\sigma_2}=\frac{1}{\sigma}+1\quad\text{and}\quad\frac{1}{\tau_1}+\frac{1}{\tau_2}\ge\frac{1}{\tau}
\end{equation}
with the convention that $1/\infty=0$. Then for any $f_1\in L^{\sigma_1,\tau_1}\(\R^n\)$ and $f_2\in L^{\sigma_2,\tau_2}\(\R^n\)$, we have $f_1\ast f_2\in L^{\sigma,\tau}\(\R^n\)$ and
$$\left\|f_1\ast f_2\right\|_{L^{\sigma,\tau}\(\R^n\)}\le C\left\|f_1\right\|_{L^{\sigma_1,\tau_1}\(\R^n\)}\left\|f_2\right\|_{L^{\sigma_2,\tau_2}\(\R^n\)}$$
for some constant $C$ independent of $\Omega$, $f_1$, and $f_2$.
\end{lemma}

\begin{lemma}{\rm(A limit case in Young's inequality)}\label{Lem3}
For any $\sigma\in\(1,\infty\)$, $f_1\in L^1\(\R^n\)$, and $f_2\in L^{\sigma,\infty}\(\R^n\)$, we have $f_1\ast f_2\in L^{\sigma,\infty}\(\R^n\)$ and
$$\left\|f_1\ast f_2\right\|_{L^{\sigma,\infty}\(\R^n\)}\le C\left\|f_1\right\|_{L^1\(\R^n\)}\left\|f_2\right\|_{L^{\sigma,\infty}\(\R^n\)}$$
for some constant $C$ independent of $f_1$ and $f_2$.
\end{lemma}

The main result of this section is the following:

\begin{lemma}\label{Lem4}
Assume that \eqref{Eq2} and \eqref{Eq5} hold true. Then for any solution $\(u,v\)\in L^a\(\R^n\)\times L^b\(\R^n\)$ of \eqref{Eq1} such that \eqref{Eq8} holds true, we have $u,v\in C^2\(\R^n\)\cap L^\infty\(\R^n\)$, $u\in L^{n/\(n-2\),\infty}\(\R^n\)$, and 
\begin{equation}\label{Lem4Eq1}
v\in\left\{\begin{aligned}&L^{n/\(n-2\),\infty}\(\R^n\)&&\text{if }q+s>n/\(n-2\)\\&L^{n\(1-s\)/\(\(n-2\)q-2\),\infty}\(\R^n\)&&\text{if }q+s<n/\(n-2\).\end{aligned}\right.
\end{equation}
In case $q+s=n/\(n-2\)$, we have $v\in L^{\sigma,\infty}\(\R^n\)$ for all $\sigma>n/\(n-2\)$ and
\begin{equation}\label{Lem4Eq2}
\left\|v\right\|_{L^{\sigma,\infty}\(\R^n\)}\le\Lambda_0\big(\(\(n-2\)\sigma-n\)^{-1/\(1-s\)}+1\big)
\end{equation}
for some constant $\Lambda_0$ independent of $\sigma$.
\end{lemma}

\proof[Proof of Lemma~\ref{Lem4}]
For this proof we borrow some ideas from Jannelli and Solimini~\cite{JanSol}. We fix a solution $\(u,v\)\in L^a\(\R^n\)\times L^b\(\R^n\)$ of \eqref{Eq1}. We let $\Gamma$ be the fundamental solution defined as
$$\Gamma\(x\):=\(\(n-2\)\omega_{n-1}\)^{-1}\left|x\right|^{2-n}\qquad\forall x\in\R^n\backslash\left\{0\right\}$$
where $\omega_n$ is the volume of the unit $\(n-1\)$--dimensional sphere. 

\smallskip
As a first step, we prove the following result:

\begin{step}\label{Lem4St1}
$u\equiv\Gamma\ast\(\left|v\right|^p\left|u\right|^{r-1}u\)$ and $v\equiv\Gamma\ast\(\left|u\right|^q\left|v\right|^{s-1}v\)$ in $\R^n.$
\end{step}

\proof[Proof of Step~\ref{Lem4St1}]
Since $\Gamma\in L^{n/\(n-2\),\infty}\(\R^n\)$, $u\in  L^a\(\R^n\)$, and $v\in L^b\(\R^n\)$, by applying Lemma~\ref{Lem2} and using \eqref{Eq3} and \eqref{Eq5}, we obtain $\widetilde{u}:=\Gamma\ast\(\left|v\right|^p\left|u\right|^{r-1}u\)\in L^a\(\R^n\)$ and $\widetilde{v}:=\Gamma\ast\(\left|u\right|^q\left|v\right|^{s-1}v\)\in L^b\(\R^n\)$. Since moreover the equations $-\Delta\widetilde{u}=\left|v\right|^p\left|u\right|^{r-1}u$ and $-\Delta\widetilde{v}=\left|u\right|^q\left|v\right|^{s-1}v$ are satisfied in the distributional sense, we obtain that the functions $u-\widetilde{u}$ and $v-\widetilde{v}$ are harmonic in $\R^n$. It follows from Greene and Wu's Liouville-type result for harmonic functions in $L^p$--spaces~\cite{GreWu} that $u\equiv\widetilde{u}$ and $v\equiv\widetilde{v}$ in $\R^n$. This ends the proof of Step~\ref{Lem4St1}.
\endproof

Then we prove the following result:

\begin{step}\label{Lem4St2}
There exists $\varepsilon>0$ such that for any $\sigma\in\(n/\(n-2\),a+\varepsilon\)$, we have $u\in L^\sigma\(\R^n\)$.
\end{step}

\proof[Proof of Step~\ref{Lem4St2}]
We fix $\varepsilon>0$ and $\sigma\in\(n/\(n-2\),a+\varepsilon\)$. We assume by contradiction that $u\not\in L^\sigma\(\R^n\)$. Integration theory yields that there exists a sequence $\(\varphi_\alpha\)_{\alpha\in\N}$ in $C^\infty_c\(\R^n\)$ such that
\begin{equation}\label{Lem4St2Eq1}
\lim_{\alpha\to\infty}\int_{\R^n}u\varphi_\alpha dx=\infty\quad\text{and}\quad\left\|\varphi_\alpha\right\|_{L^{\sigma'}\(\R^n\)}\le1\quad\forall\alpha\in\N
\end{equation}
where $\sigma':=\sigma/\(\sigma-1\)$. Moreover, since $u\in L^a\(\R^n\)$, we can choose $\varphi_\alpha$ so that 
\begin{equation}\label{Lem4St2Eq2}
2\int_{\R^n}u\varphi_\alpha dx\ge S_\alpha:=\sup_{u\in K_\alpha}\int_{\R^n}u\varphi dx
\end{equation}
where
\begin{multline*}
K_\alpha:=\bigg\{\varphi\in L^{\sigma'}\(\R^n\)\cap L^{a'}\(\R^n\);\quad\left\|\varphi\right\|_{L^{\sigma'}\(\R^n\)}\le\left\|\varphi_\alpha\right\|_{L^{\sigma'}\(\R^n\)}\\
\text{and}\quad\left\|\varphi\right\|_{L^{a'}\(\R^n\)}\le\left\|\varphi_\alpha\right\|_{L^{a'}\(\R^n\)}\bigg\}.
\end{multline*}
On the other hand, since $\(u,v\)\in L^a\(\R^n\)\times L^b\(\R^n\)$, $b\ge a$, and $p\ge1$, by applying H\"older's inequality and using \eqref{Eq3}, we obtain $\left|v\right|^{p-2}v\left|u\right|^{r-1}u\in L^\mu\(\R^n\)$ where
$$\frac{1}{\mu}:=\frac{1}{a}-\frac{1}{b}+\frac{2}{n}\,.$$
It follows that there exist sequences $\(f_\beta\)_{\beta\in\N}$ in $C^\infty_c\(\R^n\)$ and $\(g_\beta\)_{\beta\in\N}$ in $L^\mu\(\R^n\)$ such that 
\begin{equation}\label{Lem4St2Eq3}
\lim_{\beta\to\infty}\left\|g_\beta\right\|_{L^\mu\(\R^n\)}=0\quad\text{and}\quad \left|v\right|^{p-2}v\left|u\right|^{r-1}u=f_\beta+g_\beta\quad\forall\beta\in\N.
\end{equation}
It follows from \eqref{Lem4St2Eq3} and Step~\ref{Lem4St1} that
\begin{equation}\label{Lem4St2Eq4}
\int_{\R^n}u\varphi_\alpha dx=\int_{\R^n}\(\Gamma\ast\(f_\beta v\)\)\varphi_\alpha dx+\int_{\R^n}\(\Gamma\ast\(g_\beta v\)\)\varphi_\alpha dx.
\end{equation}
Since $\sigma>n/\(n-2\)$, by applying Lemmas~\ref{Lem1} and~\ref{Lem2} and using \eqref{Lem4St2Eq1}, we obtain 
\begin{align}\label{Lem4St2Eq5}
&\int_{\R^n}\(\Gamma\ast\(f_\beta v\)\)\varphi_\alpha dx\nonumber\\
&\qquad=\OO\big(\left\|\Gamma\right\|_{L^{n/\(n-2\),\infty}\(\R^n\)}\left\|f_\beta v\right\|_{L^{n\sigma/\(n+2\sigma\),\sigma}\(\R^n\)}\left\|\varphi_\alpha\right\|_{L^{\sigma'}\(\R^n\)}\big)\nonumber\\
&\qquad=\OO\big(\left\|f_\beta v\right\|_{L^{n\sigma/\(n+2\sigma\)
,\sigma}\(\R^n\)}\big).
\end{align} 
With regard to the second integral in the right-hand side of \eqref{Lem4St2Eq4}, by applying Fubini's theorem and Step~\ref{Lem4St1}, we obtain
\begin{align}\label{Lem4St2Eq6}
\int_{\R^n}\(\Gamma\ast\(g_\beta v\)\)\varphi_\alpha dx&=\int_{\R^n}g_\beta v\(\Gamma\ast\varphi_\alpha\)dx\nonumber\\
&=\int_{\R^n}g_\beta\(\Gamma\ast\(\left|u\right|^q\left|v\right|^{s-1}v\)\)\(\Gamma\ast\varphi_\alpha\)dx\nonumber\\
&=\int_{\R^n}\(\Gamma\ast\(g_\beta\(\Gamma\ast\varphi_\alpha\)\)\)\left|u\right|^q\left|v\right|^{s-1}vdx.
\end{align} 
Since $\sigma>n/\(n-2\)$, by applying Lemma~\ref{Lem2}, we obtain
\begin{align}\label{Lem4St2Eq7}
\left\|\Gamma\ast\varphi_\alpha\right\|_{L^{n\sigma/\(\(n-2\)\sigma-n\),\sigma'}\(\R^n\)}&=\OO\big(\left\|\Gamma\right\|_{L^{n/\(n-2\),\infty}\(\R^n\)}\left\|\varphi_\alpha\right\|_{L^{\sigma'}\(\R^n\)}\big)\nonumber\\
&=\OO\big(\left\|\varphi_\alpha\right\|_{L^{\sigma'}\(\R^n\)}\big).
\end{align}
By letting $\varepsilon$ be small enough so that 
$$\frac{1}{a}-\frac{1}{b}<\frac{1}{a+\varepsilon}\,,$$
applying Lemmas~\ref{Lem1} and~\ref{Lem2}, and using \eqref{Lem4St2Eq3} and \eqref{Lem4St2Eq7}, we obtain that if $\sigma\in\(n/\(n-2\),a+\varepsilon\)$, then
\begin{align}\label{Lem4St2Eq8}
&\left\|\Gamma\ast\(g_\beta\(\Gamma\ast\varphi_\alpha\)\)\right\|_{L^{\nu,\sigma'}\(\R^n\)}\nonumber\\
&\quad=\OO\big(\left\|\Gamma\right\|_{L^{n/\(n-2\),\infty}\(\R^n\)}\left\|g_\beta\right\|_{L^{\mu,\infty}\(\R^n\)}\left\|\Gamma\ast\varphi_\alpha\right\|_{L^{n\sigma/\(\(n-2\)\sigma-n\),\sigma'}\(\R^n\)}\big)\nonumber\\
&\quad=\oo\big(\left\|\varphi_\alpha\right\|_{L^{\sigma'}\(\R^n\)}\big)
\end{align}
as $\beta\to\infty$ where
$$\frac{1}{\nu}:=\frac{n-2}{n}+\frac{1}{a}-\frac{1}{b}-\frac{1}{\sigma}\,.$$
Since $q\ge1$, it follows from \eqref{Eq3}, \eqref{Lem4St2Eq8}, and another application of Lemma~\ref{Lem1} that
\begin{align}\label{Lem4St2Eq9}
&\left\|\(\Gamma\ast\(g_\beta\(\Gamma\ast\varphi_\alpha\)\)\)\left|u\right|^{q-2}u\left|v\right|^{s-1}v\right\|_{L^{\sigma'}\(\R^n\)}\nonumber\\
&\qquad=\OO\big(\left\|\(\Gamma\ast\(g_\beta\(\Gamma\ast\varphi_\alpha\)\)\)\right\|_{L^{\nu,\sigma'}\(\R^n\)}\left\|u\right\|^{q-1}_{L^{a,\infty}\(\R^n\)}\left\|v\right\|^{s}_{L^{b,\infty}\(\R^n\)}\big)\nonumber\\
&\qquad=\oo\big(\left\|\varphi_\alpha\right\|_{L^{\sigma'}\(\R^n\)}\big)
\end{align} 
as $\beta\to\infty$. Similarly as in \eqref{Lem4St2Eq7}--\eqref{Lem4St2Eq9}, we obtain
\begin{equation}\label{Lem4St2Eq10}
\left\|\(\Gamma\ast\(g_\beta\(\Gamma\ast\varphi_\alpha\)\)\)\left|u\right|^{q-2}u\left|v\right|^{s-1}v\right\|_{L^{a'}\(\R^n\)}=\oo\big(\left\|\varphi_\alpha\right\|_{L^{a'}\(\R^n\)}\big)
\end{equation} 
as $\beta\to\infty$. It follows from \eqref{Lem4St2Eq2}, \eqref{Lem4St2Eq6}, \eqref{Lem4St2Eq9}, and \eqref{Lem4St2Eq10} that
\begin{equation}\label{Lem4St2Eq11}
\int_{\R^n}\(\Gamma\ast\(g_\beta v\)\)\varphi_\alpha dx=\oo\(\int_{\R^n}u\varphi_\alpha dx\)
\end{equation}
as $\beta\to\infty$. Now by putting together \eqref{Lem4St2Eq4}, \eqref{Lem4St2Eq5}, and \eqref{Lem4St2Eq11}, we obtain
\begin{equation}\label{Lem4St2Eq12}
\int_{\R^n}u\varphi_\alpha dx=\OO\big(\left\|f_\beta v\right\|_{L^{n\sigma/\(n+2\sigma\),\sigma}\(\R^n\)}\big)
\end{equation}
for large $\beta$. Since $f_\beta\in C^\infty_c\(\R^n\)$ and $v\in L^b\(\R^n\)$, we obtain a contradiction between \eqref{Lem4St2Eq1} and \eqref{Lem4St2Eq12} provided $b>n\sigma/\(n+2\sigma\)$. Moreover, the latter inequality holds true when $\sigma<a+\varepsilon$ for small $\varepsilon>0$ as a consequence of $b\ge a$. This ends the proof of Step~\ref{Lem4St2}.
\endproof

Then we prove the following result:

\begin{step}\label{Lem4St3}
$u\in L^{n/\(n-2\),\infty}\(\R^n\)$.
\end{step}

\proof[Proof of Step~\ref{Lem4St3}]
We let $\(f_\beta\)_{\beta\in\N}$ and $\(g_\beta\)_{\beta\in\N}$ be as in \eqref{Lem4St2Eq3}. It follows from Step~\ref{Lem4St1} that 
\begin{equation}\label{Lem4St3Eq1}
u\equiv\Gamma\ast\(f_\beta v\)+\Gamma\ast\(g_\beta\(\Gamma\ast\(\left|u\right|^q\left|v\right|^{s-1}v\)\)\)\quad\text{in }\R^n.
\end{equation}
For any $\sigma\in\(n/\(n-2\),a\)$, by applying Lemmas~\ref{Lem1} and~\ref{Lem2} and Step~\ref{Lem4St2} and using \eqref{Eq5}, we obtain 
\begin{align}\label{Lem4St3Eq2}
&\left\|\Gamma\ast\(\left|u\right|^q\left|v\right|^{s-1}v\)\right\|_{L^{\omega,\infty}\(\R^n\)}\nonumber\\
&\qquad=\OO\big(\left\|\Gamma\right\|_{L^{n/\(n-2\),\infty}\(\R^n\)}\left\|u\right\|^{q-1}_{L^{a,\infty}\(\R^n\)}\left\|v\right\|^s_{L^{b,\infty}\(\R^n\)}\left\|u\right\|_{L^{\sigma,\infty}\(\R^n\)}\big)\nonumber\\
&\qquad=\OO\big(1\big)
\end{align}
and then
\begin{align}\label{Lem4St3Eq3}
&\left\|\(\Gamma\ast\(g_\beta\(\Gamma\ast\(\left|u\right|^q\left|v\right|^{s-1}v\)\)\)\)\right\|_{L^{\sigma,\infty}\(\R^n\)}\nonumber\\
&\qquad=\OO\big(\left\|\Gamma\right\|_{L^{n/\(n-2\),\infty}\(\R^n\)}\left\|g_\beta\right\|_{L^{\mu,\infty}\(\R^n\)}\left\|\Gamma\ast\(\left|u\right|^q\left|v\right|^{s-1}v\)\right\|_{L^{\omega,\infty}\(\R^n\)}\big)\nonumber\\
&\qquad=\oo\big(1\big)
\end{align}
as $\beta\to\infty$ where
$$\frac{1}{\omega}:=\frac{1}{\sigma}+\frac{1}{b}-\frac{1}{a}\,.$$
It follows from \eqref{Lem4St3Eq1}--\eqref{Lem4St3Eq3} that
\begin{equation}\label{Lem4St3Eq4}
\left\|u\right\|_{L^{\sigma,\infty}\(\R^n\)}\le2\left\|\Gamma\ast\(f_\beta v\)\right\|_{L^{\sigma,\infty}\(\R^n\)}
\end{equation}
for large $\beta$. By passing to the limit into \eqref{Lem4St3Eq4} as $\sigma\to n/\(n-2\)$, we obtain
\begin{equation}\label{Lem4St3Eq5}
\left\|u\right\|_{L^{n/\(n-2\),\infty}\(\R^n\)}\le2\left\|\Gamma\ast\(f_\beta v\)\right\|_{L^{n/\(n-2\),\infty}\(\R^n\)}.
\end{equation}
Since $\Gamma\in L^{n/\(n-2\),\infty}\(\R^n\)$, $f_\beta\in C^\infty_c\(\R^n\)$, and $v\in L^1_{\loc}\(\R^n\)$, by applying Lemma~\ref{Lem3}, we obtain $\Gamma\ast\(f_\beta v\)\in L^{n/\(n-2\),\infty}\(\R^n\)$. Hence it follows from \eqref{Lem4St3Eq5} that $u\in L^{n/\(n-2\),\infty}\(\R^n\)$. This ends the proof of Step~\ref{Lem4St3}.
\endproof

Our next step is as follows:

\begin{step}\label{Lem4St4}
Let $\widehat\sigma,\overline\sigma\in\(0,\infty\)$ be such that $u\in L^{\widehat\sigma,\infty}\(\R^n\)$ and $v\in L^{\overline\sigma,\infty}\(\R^n\)$. Then the following results hold true:
\begin{enumerate}
\item[(i)] If
\begin{equation}\label{Lem4St4Eq1}
\frac{2}{n}<\frac{p}{\overline\sigma}+\frac{r}{\widehat\sigma}<1,
\end{equation}
then
\begin{equation}\label{Lem4St4Eq2}
\left\|u\right\|_{L^{\sigma,\infty}\(\R^n\)}\le C\left\|v\right\|_{L^{\overline\sigma,\infty}\(\R^n\)}^p\left\|u\right\|_{L^{\widehat\sigma,\infty}\(\R^n\)}^r
\end{equation}
where 
$$\frac{1}{\sigma}:=\frac{p}{\overline\sigma}+\frac{r}{\widehat\sigma}-\frac{2}{n}\,.$$
\item[(ii)] If
\begin{equation}\label{Lem4St4Eq3}
\frac{2}{n}<\frac{q}{\widehat\sigma}+\frac{s}{\overline\sigma}<1,
\end{equation}
then 
\begin{equation}\label{Lem4St4Eq4}
\left\|v\right\|_{L^{\sigma,\infty}\(\R^n\)}\le C\left\|u\right\|_{L^{\widehat\sigma,\infty}\(\R^n\)}^q\left\|v\right\|_{L^{\overline\sigma,\infty}\(\R^n\)}^s
\end{equation}
where 
$$\frac{1}{\sigma}=\frac{q}{\widehat\sigma}+\frac{s}{\overline\sigma}-\frac{2}{n}\,.$$
\end{enumerate}
Moreover, the constants $C$ in \eqref{Lem4St4Eq2} and \eqref{Lem4St4Eq4} are such that
\begin{equation}\label{Lem4St4Eq5}
C\le C_n\(\frac{n\sigma}{n+2\sigma}\)'=\frac{C_nn\sigma}{\(n-2\)\sigma-n}
\end{equation}
for some constant $C_n$ depending only on $n$.
\end{step}

\proof[Proof of Step~\ref{Lem4St4}]
Since $\Gamma\in L^{n/\(n-2\),\infty}\(\R^n\)$, the estimates \eqref{Lem4St4Eq2} and \eqref{Lem4St4Eq4} follow from Step~\ref{Lem4St1} by applying Lemmas~\ref{Lem1} and~\ref{Lem2}. To obtain the estimate \eqref{Lem4St4Eq5}, we remark that 
\begin{equation}\label{Lem4St4Eq6}
\left\|f\right\|_{L^{\sigma,\infty}\(\R^n\)}\le\sup_{0<\left|E\right|<\infty}\left|E\right|^{-1/\sigma'}\int_E\left|f\right|dx\le\sigma'\left\|f\right\|_{L^{\sigma,\infty}\(\R^n\)}
\end{equation}
for all $f\in L^{\sigma,\infty}\(\R^n\)$ and $\sigma\in\(1,\infty\)$ (see for instance Grafakos~\cite{Gra}*{Exercises 1.1.12}). It follows from \eqref{Lem4St4Eq6} and the standard Young's inequality in $L^1\(\R^n\)$ that
\begin{equation}\label{Lem4St4Eq7}
\left\|f_1\ast f_2\right\|_{L^{\sigma,\infty}\(\R^n\)}\le\sigma_1'\sigma_2'\left\|f_1\right\|_{L^{\sigma_1,\infty}\(\R^n\)}\left\|f_2\right\|_{L^{\sigma_2,\infty}\(\R^n\)}
\end{equation}
for all $f_1\in L^{\sigma_1,\infty}\(\R^n\)$, $f_2\in L^{\sigma_2,\infty}\(\R^n\)$, and $\sigma,\sigma_1,\sigma_2\in\(1,\infty\)$ such that \eqref{Lem2Eq} holds true. Moreover, we have 
\begin{equation}\label{Lem4St4Eq8}
\left\|f_1f_2\right\|_{L^{\sigma,\infty}\(\R^n\)}\le\frac{\sigma_1^{1/\sigma_1}\sigma_2^{1/\sigma_2}}{\sigma^{1/\sigma}}\left\|f_1\right\|_{L^{\sigma_1,\infty}\(\R^n\)}\left\|f_2\right\|_{L^{\sigma_2,\infty}\(\R^n\)}
\end{equation}
for all $f_1\in L^{\sigma_1,\infty}\(\R^n\)$, $f_2\in L^{\sigma_2,\infty}\(\R^n\)$, and $\sigma,\sigma_1,\sigma_2\in\(0,\infty\)$ such that \eqref{Lem1Eq} holds true (see for instance Grafakos~\cite{Gra}*{Exercise 1.1.15}). The estimate \eqref{Lem4St4Eq5} then follows from \eqref{Lem4St4Eq7} and \eqref{Lem4St4Eq8}.
\endproof

Then we prove the following result:

\begin{step}\label{Lem4St5}
$u,v\in C^2\(\R^n\)\cap L^\infty\(\R^n\)$.
\end{step}

\proof[Proof of Step~\ref{Lem4St5}]
It follows from Step~\ref{Lem4St2} that $u\in L^{a+\varepsilon}\(\R^n\)$ for small $\epsilon>0$. We let $\(\widehat{\sigma}_\alpha\)_{\alpha\in\N}$ and $\(\overline{\sigma}_\alpha\)_{\alpha\in\N}$ be the sequences defined by $\widehat{\sigma}_0:=a+\varepsilon$, $\overline{\sigma}_0:=b$, and
$$\left\{\begin{aligned}&\frac{1}{\overline{\sigma}_{\alpha+1}}:=\max\(\frac{q}{\widehat{\sigma}_\alpha}+\frac{s}{\overline{\sigma}_\alpha}-\frac{2}{n},0\)\\&\frac{1}{\widehat{\sigma}_{\alpha+1}}:=\max\(\frac{p}{\overline{\sigma}_{\alpha+1}}+\frac{r}{\widehat{\sigma}_\alpha}-\frac{2}{n},0\)\end{aligned}\right.\qquad\forall\alpha\in\N$$
with the convention that $1/\infty=0$. By using \eqref{Eq2} and \eqref{Eq3}, one can see that $\widehat{\sigma}_\alpha,\overline{\sigma}_\alpha\nearrow\infty$ as $\alpha\to\infty$. Moreover, by applying Step~\ref{Lem4St4} together with an induction argument, we obtain that $\(u,v\)\in L^{\widehat{\sigma},\infty}\(\R^n\)\times L^{\overline{\sigma},\infty}\(\R^n\)$ for all $\widehat{\sigma}\in\[a,\widehat{\sigma}_{\alpha}\)$ and $\overline{\sigma}\in\[b,\overline{\sigma}_{\alpha}\)$. Hence we obtain that $u,v\in L^{\sigma}\(\R^n\)$ for large $\sigma>0$. By applying Step~\ref{Lem4St1}, Calderon--Zygmund's inequality, and an interpolation inequality (see Gilbarg--Trudinger~\cite{GilTru}*{Theorems~7.28 and~9.9}), it follows that $u,v\in W^{2,\sigma}\(\R^n\)$ for large $\sigma>0$. Finally, by applying Sobolev's inequalities and standard elliptic regularity theory, we obtain $u,v\in L^\infty\(\R^n\)\cap C^2\(\R^n\)$. 
\endproof

Now we can prove the following result:

\begin{step}\label{Lem4St6}
\eqref{Lem4Eq1} and \eqref{Lem4Eq2} hold true.
\end{step}

\proof[Proof of Step~\ref{Lem4St6}] 
We separate the following cases: $q+s<n/\(n-2\)$, $q+s=n/\(n-2\)$, [$q+s>n/\(n-2\)$, $q<n/\(n-2\)$, and $s\le1$], [$q<n/\(n-2\)$ and $s>1$], and [$q\ge n/\(n-2\)$ and $s\ne0$].

\proof[Case $q+s<n/\(n-2\)$] 
In this case, we fix $\widehat\sigma:=n/\(n-2\)$ so that \eqref{Lem4St4Eq3} holds true for all $\overline\sigma$ such that
\begin{equation}\label{Lem4St6Eq1}
\overline\sigma>\frac{ns}{n-\(n-2\)q}\,.
\end{equation}
We let $\(\sigma_\alpha\)_{\alpha\in\N}$ be the sequence defined by $\sigma_0:=b$ and
$$\frac{1}{\sigma_{\alpha+1}}:=\frac{s}{\sigma_\alpha}+\frac{\(n-2\)q-2}{n}\qquad\forall\alpha\in\N.$$
Since $q+s<n/\(n-2\)$, we obtain
\begin{equation}\label{Lem4St6Eq2}
\sigma_\alpha\searrow\frac{n\(1-s\)}{\(n-2\)q-2}>\max\(\frac{ns}{n-\(n-2\)q},\frac{n}{n-2}\)
\end{equation}
as $\alpha\to\infty$. It follows from \eqref{Lem4St6Eq2} that \eqref{Lem4St6Eq1} holds true with $\overline\sigma=\sigma_\alpha$ for all $\alpha\in\N$. By applying Step~\ref{Lem4St4}, we then obtain
\begin{equation}\label{Lem4St6Eq3}
\left\|v\right\|_{L^{\sigma_{\alpha},\infty}\(\R^n\)}\le\frac{C_nn\sigma_{\alpha}}{\(n-2\)\sigma_{\alpha}-n}\left\|u\right\|_{L^{n/\(n-2\),\infty}\(\R^n\)}^q\left\|v\right\|_{L^{\sigma_{\alpha-1},\infty}\(\R^n\)}^s
\end{equation}
for all $\alpha\ge1$. It follows from \eqref{Lem4St6Eq2} and Step~\ref{Lem4St3} that
\begin{equation}\label{Lem4St6Eq4}
\frac{C_nn\sigma_{\alpha}}{\(n-2\)\sigma_{\alpha}-n}\left\|u\right\|_{L^{n/\(n-2\),\infty}\(\R^n\)}^q\le C
\end{equation}
for some constant $C$ independent of $\alpha$. It follows from \eqref{Lem4St6Eq3} and \eqref{Lem4St6Eq4} that
\begin{equation}\label{Lem4St6Eq5}
\left\|v\right\|_{L^{\sigma_\alpha,\infty}\(\R^n\)}\le C^{1+s+\dotsb+s^{\alpha-1}}\left\|v\right\|_{L^{b,\infty}\(\R^n\)}^{s^\alpha}.
\end{equation}
By passing to the limit into \eqref{Lem4St6Eq3} as $\alpha\to\infty$, we obtain
$$\left\|v\right\|_{L^{n\(1-s\)/\(\(n-2\)q-2\),\infty}\(\R^n\)}\le C^{1/\(1-s\)},$$
namely \eqref{Lem4Eq1} holds true in case $q+s<n/\(n-2\)$.

\proof[Case $q+s=n/\(n-2\)$] 
In this case, we let $\widehat\sigma$ and $\(\sigma_\alpha\)_{\alpha\in\N}$ be as in the previous case but instead of \eqref{Lem4St6Eq2}, we obtain 
\begin{equation}\label{Lem4St6Eq6}
\sigma_\alpha\searrow\frac{n\(1-s\)}{\(n-2\)q-2}=\frac{ns}{n-\(n-2\)q}=\frac{n}{n-2}
\end{equation}
as $\alpha\to\infty$. Moreover, one can check that
\begin{align}
\(n-2\)\sigma_\alpha-n&=\frac{n\(\(n-2\)b-n\)s^\alpha}{\(n-2\)b-\(\(n-2\)b-n\)s^\alpha}\nonumber\\
&\sim\frac{n\(\(n-2\)b-n\)}{\(n-2\)b}\,s^\alpha\label{Lem4St6Eq7}
\end{align}
as $\alpha\to\infty$. It follows from \eqref{Lem4St6Eq6}, \eqref{Lem4St6Eq7}, and Step~\ref{Lem4St3} that
\begin{equation}\label{Lem4St6Eq8}
\frac{C_nn\sigma_{\alpha}\left\|u\right\|_{L^{n/\(n-2\),\infty}\(\R^n\)}^q}{\(n-2\)\sigma_{\alpha}-n}\le Cs^{-\alpha}
\end{equation}
for some constant $C$ independent of $\alpha$. It follows from \eqref{Lem4St6Eq3} and \eqref{Lem4St6Eq8} that
\begin{align}
\left\|v\right\|_{L^{\sigma_\alpha,\infty}\(\R^n\)}&\le C^{1+s+\dotsb+s^{\alpha-1}}s^{-\alpha-\(\alpha-1\)s-\dotsb-s^{\alpha-1}}\left\|v\right\|_{L^{b,\infty}\(\R^n\)}^{s^\alpha}\nonumber\\
&\sim C^{1/\(1-s\)}s^{\(s-\(1-s\)\alpha\)/\(1-s\)^2}\nonumber\\
&\sim C'\(\(n-2\)\sigma_\alpha-n\)^{-1/\(1-s\)}\label{Lem4St6Eq9}
\end{align}
for some constant $C'$ independent of $\alpha$. The estimate \eqref{Lem4Eq2} then follows from \eqref{Lem4St6Eq9} by an easy interpolation argument.

\proof[Case $q+s>n/\(n-2\)$, $q<n/\(n-2\)$, and $s\le1$] 
In this case, we \linebreak choose $\widehat\sigma$ so that
\begin{equation}\label{Lem4St6Eq10}
q<\widehat\sigma<\frac{nq}{n-\(n-2\)s}\,.
\end{equation}
It follows from \eqref{Lem4St6Eq10} that 
\begin{equation}\label{Lem4St6Eq11}
\frac{s\widehat\sigma}{\widehat\sigma-q}<\frac{ns}{n-\(n-2\)q}
\end{equation}
and \eqref{Lem4St4Eq3} holds true for all $\overline\sigma$ such that
\begin{equation}\label{Lem4St6Eq12}
\overline\sigma>\frac{s\widehat\sigma}{\widehat\sigma-q}\,.
\end{equation}
Moreover, $u\in L^{\widehat\sigma,\infty}\(\R^n\)$ as a consequence of \eqref{Lem4St6Eq10} and Steps~\ref{Lem4St3} and~\ref{Lem4St5}. We let $\(\sigma_\alpha\)_{\alpha\in\N}$ be the sequence defined by
$$\frac{1}{\sigma_{\alpha+1}}:=\frac{s}{\sigma_\alpha}+\frac{q}{\widehat\sigma}-\frac{2}{n}\qquad\forall\alpha\in\N$$
with $\sigma_0$ chosen large enough so that $v\in L^{\sigma_0,\infty}\(\R^n\)$. It follows from the second inequality in \eqref{Lem4St6Eq10} that 
\begin{equation}\label{Lem4St6Eq13}
\sigma_\alpha\searrow\frac{n\widehat\sigma\(1-s\)}{nq-2\widehat\sigma}<\frac{n}{n-2}
\end{equation}
as $\alpha\to\infty$. Moreover, it follows from \eqref{Lem4St6Eq11} that we can choose $\sigma_0$ so that 
\begin{equation}\label{Lem4St6Eq14}
\frac{s\widehat\sigma}{\widehat\sigma-q}<\sigma_{\alpha_0}<\frac{ns}{n-\(n-2\)q}
\end{equation}
for some $\alpha_0\ge1$. It follows from \eqref{Lem4St6Eq13} and \eqref{Lem4St6Eq14} that \eqref{Lem4St6Eq12} holds true with $\overline\sigma=\sigma_\alpha$ for all $\alpha\in\left\{0,\dotsc,\alpha_0\right\}$. By applying Step~\ref{Lem4St4}, we then obtain
\begin{equation}\label{Lem4St6Eq15}
\left\|v\right\|_{L^{\sigma_{\alpha},\infty}\(\R^n\)}=\OO\big(\left\|u\right\|_{L^{\widehat\sigma,\infty}\(\R^n\)}^q\left\|v\right\|_{L^{\sigma_{\alpha-1},\infty}\(\R^n\)}^s\big)
\end{equation}
for all $\alpha\in\left\{1,\dotsc,\alpha_0\right\}$. It follows from \eqref{Lem4St6Eq15} that
\begin{equation}\label{Lem4St6Eq16}
\left\|v\right\|_{L^{\sigma_{\alpha_0},\infty}\(\R^n\)}=\OO\big(\left\|u\right\|_{L^{\widehat\sigma,\infty}\(\R^n\)}^{q(1+s+\dotsb+s^{\alpha_0-1})}\left\|v\right\|_{L^{\sigma_0,\infty}\(\R^n\)}^{s^{\alpha_0}}\big)=\OO\(1\).
\end{equation}
It follows from \eqref{Lem4St6Eq14} that we can choose $\overline\sigma_0$ so that
\begin{equation}\label{Lem4St6Eq17}
\sigma_{\alpha_0}<\overline\sigma_0<\frac{ns}{n-\(n-2\)q}
\end{equation}
By using \eqref{Lem4St6Eq16}, the first inequality in \eqref{Lem4St6Eq17}, and Step~\ref{Lem4St5}, we obtain $v\in L^{\overline\sigma_0}\(\R^n\)$. Moreover, by using the second inequality in \eqref{Lem4St6Eq17}, we obtain
\begin{equation}\label{Lem4St6Eq18}
\widehat\sigma_0:=\frac{q\overline\sigma_0}{\overline\sigma_0-s}>\frac{n}{n-2}\,.
\end{equation}
It follows from \eqref{Lem4St6Eq18} and Steps \ref{Lem4St3} and~\ref{Lem4St5} that $u\in L^{\widehat\sigma_0}\(\R^n\)$ . By applying Lemmas~\ref{Lem1} and~\ref{Lem3} and Step~\ref{Lem4St1}, we then obtain
$$\left\|v\right\|_{L^{n/\(n-2\),\infty}\(\R^n\)}=\OO\big(\left\|\Gamma\right\|_{L^{n/\(n-2\),\infty}\(\R^n\)}\left\|u\right\|_{L^{\widehat\sigma_0}\(\R^n\)}^q\left\|v\right\|_{L^{\overline\sigma_0}\(\R^n\)}^s\big)=\OO\(1\),$$
namely \eqref{Lem4Eq1} holds true in case $q+s>n/\(n-2\)$, $q<n/\(n-2\)$, and $s\le1$.

\proof[Case $q<n/\(n-2\)$ and $s>1$]
Note that by using \eqref{Eq2}--\eqref{Eq5}, we obtain that $s>1$ implies that $q+s>n/\(n-2\)$. In this case, we choose $\widehat\sigma$ so that
\begin{equation}\label{Lem4St6Eq19}
\frac{1}{a}<\frac{1}{\widehat\sigma}<\min\(\frac{n-2}{n},\frac{1}{a}+\frac{n-2}{nq},\frac{1}{a}+\frac{\(n-2\)\(q+s\)-n}{nsq}\)\,.
\end{equation}
By using \eqref{Lem4St6Eq19} and Step~\ref{Lem4St3} and since $u\in L^a\(\R^n\)$, we obtain $u\in L^{\widehat\sigma,\infty}\(\R^n\)$. We let $\(\sigma_\alpha\)_{\alpha\in\N}$ be the sequence defined by
$$\frac{1}{\sigma_{\alpha+1}}:=\frac{1}{\sigma_\alpha}+\frac{s-1}{b}+\frac{q}{\widehat\sigma}-\frac{2}{n}\qquad\forall\alpha\in\N$$
with $\sigma_0$ chosen large enough so that $v\in L^{\sigma_0,\infty}\(\R^n\)$. It follows from \eqref{Eq3} and \eqref{Lem4St6Eq19} that $\sigma_\alpha\searrow0$ as $\alpha\to\infty$ and that we can choose $\sigma_0$ so that 
\begin{equation}\label{Lem4St6Eq20}
\frac{n-\(n-2\)q}{ns}<\frac{1}{\sigma_{\alpha_0}}<\frac{q}{a}-\frac{q}{\widehat\sigma}+\frac{n-2}{n}
\end{equation}
for some $\alpha_0\ge1$. It follows from \eqref{Eq3}, \eqref{Lem4St6Eq19}, and \eqref{Lem4St6Eq20} that \eqref{Lem4St4Eq3} holds true with $\overline\sigma=sb\sigma_\alpha/\(\(s-1\)\sigma_\alpha+b\)$ for all $\alpha\in\left\{0,\dotsc,\alpha_0\right\}$. By applying Step~\ref{Lem4St4} and Lemma~\ref{Lem1}, we then obtain
\begin{align}
\left\|v\right\|_{L^{\sigma_\alpha,\infty}\(\R^n\)}&=\OO\big(\left\|u\right\|_{L^{\widehat\sigma,\infty}\(\R^n\)}^q\left\|v\right\|_{L^{sb\sigma_{\alpha-1}/\(\(s-1\)\sigma_{\alpha-1}+b\),\infty}\(\R^n\)}^{s}\big)\nonumber\\
&=\OO\big(\left\|u\right\|_{L^{\widehat\sigma,\infty}\(\R^n\)}^q\left\|v\right\|_{L^{b,\infty}\(\R^n\)}^{s-1}\left\|v\right\|_{L^{\sigma_{\alpha-1},\infty}\(\R^n\)}\big)\label{Lem4St6Eq21}
\end{align}
for all $\alpha\in\left\{1,\dotsc,\alpha_0\right\}$. It follows from \eqref{Lem4St6Eq21} that
$$\left\|v\right\|_{L^{\sigma_{\alpha_0},\infty}\(\R^n\)}=\OO\big(\left\|u\right\|_{L^{\widehat\sigma,\infty}\(\R^n\)}^{\alpha_0q}\left\|v\right\|_{L^{b,\infty}\(\R^n\)}^{\alpha_0\(s-1\)}\left\|v\right\|_{L^{\sigma_0,\infty}\(\R^n\)}\big)\\
=\OO\(1\).$$
We then choose $\overline\sigma_0$ and $\widehat\sigma_0$ as in \eqref{Lem4St6Eq17} and \eqref{Lem4St6Eq17} and we conclude in the same way as in the previous case. This ends the proof of Step~\ref{Lem4St6} and also Lemma~\ref{Lem4}.

\proof[Case $q\ge n/\(n-2\)$ and $s\ne0$]
 In this case, we choose $\overline\sigma_0$ large enough so that $v\in L^{\overline\sigma_0}\(\R^n\)$ and we let $\widehat\sigma_0$ be as in \eqref{Lem4St6Eq17}. We then obtain $\widehat\sigma_0>n/\(n-2\)$ hence $u\in L^{\widehat\sigma_0}\(\R^n\)$. We then conclude in the same way as in the previous two cases.
\endproof

\section{Decay estimates}\label{Sec3}

In this section we establish the a priori estimates \eqref{Th1Eq1} and \eqref{Th1Eq3}. We assume that \eqref{Eq2} and \eqref{Eq5} hold true. We fix a solution $\(u,v\)\in L^a\(\R^n\)\times L^b\(\R^n\)$ of \eqref{Eq1} such that \eqref{Eq8} and \eqref{Eq9} hold true.

\smallskip
First we establish the following preliminary (non-sharp) estimate:

\begin{step}\label{Th1St1}
There exists a constant $K_0$ such that
\begin{equation}\label{Th1St1Eq1}
\left|u\(x\)\right|^a+\left|v\(x\)\right|^b\le K_0\(1+\left|x\right|^n\)^{-1}\qquad\forall x\in\R^n.
\end{equation}
\end{step}

\proof[Proof of Step~\ref{Th1St1}]
We assume by contradiction that \eqref{Th1St1Eq1} does not hold true. Since $u,v\in L^\infty\(\R^n\)$, it follows there exists a sequence of points $\(x_\alpha\)_{\alpha\in\N}$ in $\R^n$ such that 
$$\left|x_\alpha\right|M\(x_\alpha\)>2\alpha\qquad\forall\alpha\in\N$$
where
$$M\(x\):=\big(\left|u\(x\)\right|^{a}+\left|v\(x\)\right|^{b}\big)^{1/n}.$$
By applying a doubling property (see Pol\'a\v{c}ik, Quittner, and Souplet~\cite{PolQuiSou}*{Lemma~5.1}), we obtain that for any $\alpha$, there exists a point $y_\alpha\in\R^n$ such that
\begin{equation}\label{Th1St1Eq2}
\left|y_\alpha\right|M\(y_\alpha\)>2\alpha\,,\quad M\(x_\alpha\)\le M\(y_\alpha\),
\end{equation}
and 
\begin{equation}\label{Th1St1Eq3}
M\(y\)\le2M\(y_\alpha\)\qquad\forall y\in B\big(y_\alpha,\alpha/M\(y_\alpha\)\big).
\end{equation}
We define
$$\left\{\begin{aligned}u_\alpha\(y\):=M\(y_\alpha\)^{-n/a}u\big(M\(y_\alpha\)^{-1}y+y_\alpha\big)\\v_\alpha\(y\):=M\(y_\alpha\)^{-n/b}v\big(M\(y_\alpha\)^{-1}y+y_\alpha\big)\end{aligned}\right.\qquad\forall y\in\R^n.$$
By using \eqref{Eq1}, we obtain
\begin{equation}\label{Th1St1Eq4}
\left\{\begin{aligned}&-\Delta u_\alpha=\left|v_\alpha\right|^{p}\left|u_\alpha\right|^{r-1}u_\alpha&&\text{in }\R^n\\
&-\Delta v_\alpha=\left|u_\alpha\right|^{q}\left|v_\alpha\right|^{s-1}v_\alpha&&\text{in }\R^n.
\end{aligned}\right.
\end{equation}
Moreover, it follows from \eqref{Th1St1Eq3} that
\begin{equation}\label{Th1St1Eq5}
\left|u_\alpha\(y\)\right|^{a}+\left|v_\alpha\(y\)\right|^{b}\le2^n\qquad\forall y\in B\big(0,\alpha\big).
\end{equation}
By applying standard elliptic estimates, it follows from \eqref{Th1St1Eq4} and \eqref{Th1St1Eq5} that the sequences $\(u_\alpha\)_{\alpha\in\N}$ and $\(v_\alpha\)_{\alpha\in\N}$ converge up to a subsequence in $C^2_{\loc}\(\R^n\)$ to some functions $u_\infty$ and $v_\infty$, respectively. On the other hand, by using the definition of $u_\alpha$ and $v_\alpha$, we obtain
\begin{equation}\label{Th1St1Eq6}
\int_{B\(0,R\)}\big(\left|u_\alpha\right|^a+\left|v_\alpha\right|^b\big)dx=\int_{B(y_\alpha,R/M\(y_\alpha\))}\big(\left|u\right|^a+\left|v\right|^b\big)dx.
\end{equation}
for all $R>0$. Since $\(u,v\)\in L^a\(\R^n\)\times L^b\(\R^n\)$ and $u,v\in L^\infty\(\R^n\)$, it follows from the first inequality in \eqref{Th1St1Eq2} that
\begin{equation}\label{Th1St1Eq7}
\int_{B(y_\alpha,R/M\(y_\alpha\))}\big(\left|u\right|^a+\left|v\right|^b\big)dx\longrightarrow0
\end{equation}
as $\alpha\to\infty$. It follows from \eqref{Th1St1Eq6} and \eqref{Th1St1Eq7} that $u_\infty\equiv v_\infty\equiv0$ in $\R^n$. However by using the definitions of $u_\alpha$ and $v_\alpha$, we obtain 
$$\left|u_\alpha\(0\)\right|^{a}+\left|v_\alpha\(0\)\right|^{b}=1$$
which gives 
$$\left|u_\infty\(0\)\right|^{a}+\left|v_\infty\(0\)\right|^{b}=1.$$
Thus we obtain a contradiction. This ends the proof of Step~\ref{Th1St1}.
\endproof

Then we prove the following result:

\begin{step}\label{Th1St2}
If $p\ge s$ (this holds true in case $r<1$ as a consequence of \eqref{Eq2}) and $u\ge0$ in $\R^n$, then
\begin{equation}\label{Th1St2Eq1}
\frac{\left|v\(x\)\right|^{p-s+1}}{p-s+1}\le\frac{u\(x\)^{q-r+1}}{q-r+1}\qquad\forall x\in\R^n.
\end{equation}
\end{step}

\proof[Proof of Step~\ref{Th1St2}]
The inequality \eqref{Th1St2Eq1} has been obtained by Quittner and Souplet~\cite{QuiSou} in case $p+r=q+s$ and $u,v\ge0$ in $\R^n$ and Souplet~\cite{Soup} in case $r=s=0$ and $u,v\ge0$ in $\R^n$. Here we adapt the proofs of~\cite{QuiSou} and~\cite{Soup} to our setting. We assume that $p\ge s$ and $u\ge0$ in $\R^n$. By applying the strong maximum principle and since $u,v\in L^\infty\(\R^n\)$, we obtain that either $u\equiv v\equiv0$ or $u>0$ in $\R^n$. Hence we can assume that $u>0$. We define 
$$w:=\left|v\right|-cu^\theta$$
where
$$\theta:=\frac{q-r+1}{p-s+1}\quad\text{and}\quad c:=\theta^{-1/\(p-s+1\)}.$$
It follows from \eqref{Eq2} that $0<\theta\le1$. Note that since $u>0$ in $\R^n$, we obtain $\left|v\right|>0$ in $W:=\left\{w\ge0\right\}$ hence $w$ is smooth in $W$. Moreover, by using \eqref{Eq1} and since $0<\theta\le1$ and $p\ge s$, we obtain
\begin{multline}\label{Th1St2Eq2}
\Delta w=\Delta\left|v\right|+c\theta\(1-\theta\)u^{\theta-2}\left|\nabla u\right|^2-c\theta u^{\theta-1}\Delta u\ge c\theta u^{\theta+r-1}\left|v\right|^p-u^q\left|v\right|^s\\
=c\theta\left|v\right|^su^{\theta+r-1}\(\left|v\right|^{p-s}-c^{p-s}u^{\theta\(p-s\)}\)\ge0\quad\text{in }W.
\end{multline}
For any $\varepsilon>0$, we define $w_\varepsilon:=\varepsilon H\(w/\varepsilon\)$ where $H\in C^\infty\(\R\)$ is such that $H\(s\)=0$ for all $s\in\(-\infty,0\)$ and $H\(s\)=s$ for all $s\in\(1,\infty\)$. For any $\varepsilon,R>0$, by integrating by parts, we obtain
$$\int_{B\(0,R\)}\left|\nabla w_\varepsilon\right|^2dx=\int_{\partial B\(0,R\)}w_\varepsilon\frac{\partial w_\varepsilon}{\partial\nu}d\sigma-\int_{B\(0,R\)}w_\varepsilon\Delta w_\varepsilon dx$$
which gives
\begin{multline}\label{Th1St2Eq3}
\int_{B\(0,R\)}H'\(w/\varepsilon\)^2\left|\nabla w\right|^2dx=\varepsilon\int_{\partial B\(0,R\)}H\(w/\varepsilon\)H'\(w/\varepsilon\)\frac{\partial w}{\partial\nu}d\sigma\\
-\int_{B\(0,R\)}H\(w/\varepsilon\)\(\varepsilon H'\(w/\varepsilon\)\Delta w+H''\(w/\varepsilon\)\left|\nabla w\right|^2\)dx.
\end{multline}
By using \eqref{Th1St2Eq2} and passing to the limit into \eqref{Th1St2Eq3} as $\varepsilon\to0$, we obtain
\begin{align}\label{Th1St2Eq4}
\int_{B\(0,R\)\cap W}\left|\nabla w\right|^2dx&\le\int_{\partial B\(0,R\)\cap W}w\frac{\partial w}{\partial\nu}d\sigma\nonumber\\
&=\frac{1}{2}R^{n-1}\frac{d}{dr}\[\int_{\partial B\(0,1\)}w_+\(rz\)^2d\sigma\(z\)\]_{r=R}\hspace{-2pt}
\end{align}
where $w_+:=\max\(w,0\)$ and $d\sigma$ is the volume element on $\partial B\(0,1\)$ and $\partial B\(0,R\)$. On the other hand, by applying Step~\ref{Th1St1}, we obtain
\begin{equation}\label{Th1St2Eq5}
w\(x\)<\left|v\(x\)\right|\le K_0^{1/b}\(1+\left|x\right|^n\)^{-1/b}\qquad\forall x\in\R^n.
\end{equation}
It follows from \eqref{Th1St2Eq5} that there exists a sequence of positive real numbers $\(R_\alpha\)_{\alpha\in\N}$ such that $R_\alpha\to\infty$ as $\alpha\to\infty$ and 
$$\frac{d}{dr}\[\int_{\partial B\(0,1\)}w_+\(rz\)^2d\sigma\(z\)\]_{r=R_\alpha}<0.$$
By applying \eqref{Th1St2Eq4} with $R=R_\alpha$ and letting $\alpha\to\infty$, we obtain that $\nabla w\equiv0$ in $W$. It follows from \eqref{Th1St2Eq5} that $w\le0$ in $\R^n$, namely \eqref{Th1St2Eq1} holds true. This ends the proof of Step~\ref{Th1St2}.
\endproof

Now we can establish our sharp upper bound estimates:

\begin{step}\label{Th1St3}
The upper bound estimates \eqref{Th1Eq1} hold true.
\end{step}

\proof[Proof of Step~\ref{Th1St3}] 
We separate the cases $s<1$ and $s\ge1$. 

\proof[Case $s<1$] In this case, it follows from \eqref{Eq9} that either $u\ge0$ in $\R^n$ or $r\ge1$. We begin with proving the upper bound estimate for $u$. For any $y\in\R^n$ and $R\ge1$, we define 
\begin{equation}\label{Th1St3Eq1}
u_R\(y\):=R^{n-2}u\(Ry\).
\end{equation}
It follows from the first equation of \eqref{Eq1} that
\begin{equation}\label{Th1St3Eq2}
-\Delta u_R\(y\)=f_R\(y\):=R^2\left|v\(Ry\)\right|^{p}\left|u\(Ry\)\right|^{r-1}u_R\(y\).
\end{equation}
By applying Step~\ref{Th1St1}, together with Step~\ref{Th1St2} in case $r<1$, and using \eqref{Eq3}, we obtain 
\begin{equation}\label{Th1St3Eq3}
\left|f_R\(y\)\right|\le C\left|y\right|^{-2}\left|u_R\(y\)\right|
\end{equation}
for some constant $C$ independent of $y$ and $R$. Moreover, by applying Kato's inequality~\cite{Kato}, we obtain
\begin{equation}\label{Th1St3Eq4}
-\Delta \left|u_R\right|\le\left|\Delta u_R\right|\quad\text{in }\R^n
\end{equation}
where the inequality must be understood in the weak sense. By applying a weak Harnack-type inequality (see Trudinger~\cite{Tru}*{Theorem~1.3}),
it follows from \eqref{Th1St3Eq2}--\eqref{Th1St3Eq4} that for any $\sigma>1$, there exists a constant $C_\sigma$ independent of $R$ such that
\begin{equation}\label{Th1St3Eq5}
\left\|u_R\right\|_{L^\infty\(B\(0,4\)\backslash B\(0,2\)\)}\le C_\sigma\left\|u_R\right\|_{L^\sigma\(B\(0,5\)\backslash B\(0,1\)\)}.
\end{equation}
By taking $\sigma\in\(1,n/\(n-2\)\)$ and applying Lemma~\ref{Lem1}, we obtain
\begin{equation}\label{Th1St3Eq6}
\left\|u_R\right\|_{L^\sigma\(B\(0,5\)\backslash B\(0,1\)\)}\le C'_\sigma\left\|u_R\right\|_{L^{n/\(n-2\),\infty}\(B\(0,5\)\backslash B\(0,1\)\)}
\end{equation}
for some constant $C'_\sigma$ independent of $R$. We observe that 
\begin{equation}\label{Th1St3Eq7}
\left\|u_R\right\|_{L^{n/\(n-2\),\infty}\(B\(0,5\)\backslash B\(0,1\)\)}=\left\|u\right\|_{L^{n/\(n-2\),\infty}\(B\(0,5R\)\backslash B\(0,R\)\)}.
\end{equation}
It follows from \eqref{Th1St3Eq5}--\eqref{Th1St3Eq7}, and Lemma~\ref{Lem4} that
\begin{equation}\label{Th1St3Eq8}
\left\|u_R\right\|_{L^{\infty}\(B\(0,4\)\backslash B\(0,2\)\)}\le C
\end{equation} 
for some constant $C$ independent of $R$. By applying \eqref{Th1St3Eq8} with $R=\left|x\right|/3$, we obtain
\begin{equation}\label{Th1St3Eq9}
\left|u\(x\)\right|\le C\(\left|x\right|/3\)^{2-n}\qquad\forall x\in\R^n\backslash B\(0,3\).
\end{equation}
Since on the other hand, $u$ is continuous in $\R^n$, we can deduce from \eqref{Th1St3Eq9} that the upper bound estimate for $u$ in \eqref{Th1Eq1} holds true. 

\smallskip
Now we establish the upper bound estimate for $v$. For any $y\in\R^n$ and $R\ge1$, we define 
\begin{equation}\label{Th1St3Eq10}
v_R\(y\):=\left\{\begin{aligned}&R^{n-2}v\(Ry\)&&\text{if }q+s>n/\(n-2\)\\&\frac{R^{n-2}v\(Ry\)}{\ln\(1+R\)^{1/\(1-s\)}}&&\text{if }q+s=n/\(n-2\)\\&R^{\(\(n-2\)q-2\)/\(1-s\)}v\(Ry\)&&\text{if }q+s<n/\(n-2\).\end{aligned}\right.
\end{equation}
It follows from the second equation of \eqref{Eq1} that
\begin{equation}\label{Th1St3Eq11}
-\Delta v_R\(y\)=g_R\(y\):=R^2\left|u\(Ry\)\right|^{q}\left|v\(Ry\)\right|^{s-1}v_R\(y\).
\end{equation}
By using the upper bound estimate for $u$, we obtain 
\begin{align}
\left|g_R\(y\)\right|&\le\left\{\begin{aligned}&CR^{n-\(q+s\)\(n-2\)}\left|y\right|^{-q\(n-2\)}\left|v_R\(y\)\right|^s&&\text{if }q+s>n/\(n-2\)\\
&C\ln\(1+R\)^{-1}\left|y\right|^{-q\(n-2\)}\left|v_R\(y\)\right|^s&&\text{if }q+s=n/\(n-2\)\\
&C\left|y\right|^{-q\(n-2\)}\left|v_R\(y\)\right|^s&&\text{if }q+s<n/\(n-2\)\end{aligned}\right.\nonumber\\
&\le C'\left|y\right|^{-q\(n-2\)}\(\left|v_R\(y\)\right|+1\)\label{Th1St3Eq12}
\end{align}
for some constants $C,C'$ independent of $y$ and $R$. Moreover, by applying Kato's inequality~\cite{Kato}, we obtain
\begin{equation}\label{Th1St3Eq13}
-\Delta \left|v_R\right|\le\left|\Delta v_R\right|\quad\text{in }\R^n
\end{equation}
where the inequality must be understood in the weak sense. Now we separate two cases. First we consider the case $q+s\ne n/\(n-2\)$. By applying a  weak Harnack-type inequality (see Trudinger~\cite{Tru}*{Corollary~1.1}), it follows from \eqref{Th1St3Eq11}--\eqref{Th1St3Eq13} that for any $\sigma>1$, there exists a constant $C_\sigma$ independent of $R$ such that
\begin{equation}\label{Th1St3Eq14}
\left\|v_R\right\|_{L^\infty\(B\(0,4\)\backslash B\(0,2\)\)}\le C_\sigma\big(\left\|v_R\right\|_{L^\sigma\(B\(0,5\)\backslash B\(0,1\)\)}+1\big).
\end{equation}
By taking $\sigma\in\(1,\sigma_0\)$ where
$$\sigma_0:=\left\{\begin{aligned}&n/\(n-2\)&&\text{if }q+s>n/\(n-2\)\\&n\(1-s\)/\(\(n-2\)q-2\)&&\text{if }q+s<n/\(n-2\)
\end{aligned}\right.$$ 
and applying Lemma~\ref{Lem1}, we obtain
\begin{equation}\label{Th1St3Eq15}
\left\|v_R\right\|_{L^\sigma\(B\(0,5\)\backslash B\(0,1\)\)}\le C'_\sigma\left\|v_R\right\|_{L^{\sigma_0,\infty}\(B\(0,5\)\backslash B\(0,1\)\)}
\end{equation} 
for some constant $C'_\sigma$ independent of $R$. We observe that 
\begin{equation}\label{Th1St3Eq16}
\left\|v_R\right\|_{L^{\sigma_0,\infty}\(B\(0,5\)\backslash B\(0,1\)\)}=\left\|v\right\|_{L^{\sigma_0,\infty}\(B\(0,5R\)\backslash B\(0,R\)\)}.
\end{equation}
It follows from \eqref{Th1St3Eq14}--\eqref{Th1St3Eq16}, and Lemma~\ref{Lem4} that
\begin{equation}\label{Th1St3Eq17}
\left\|v_R\right\|_{L^{\infty}\(B\(0,4\)\backslash B\(0,2\)\)}\le C
\end{equation} 
for some constant $C$ independent of $R$. By applying \eqref{Th1St3Eq17} with $R=\left|x\right|/3$, we obtain
\begin{equation}\label{Th1St3Eq18}
\left|v\(x\)\right|\le C\(\left|x\right|/3\)^{-n/\sigma_0}\qquad\forall x\in\R^n\backslash B\(0,3\).
\end{equation}
Since on the other hand, $v$ is continuous in $\R^n$, we can deduce from \eqref{Th1St3Eq18} that the upper bound estimate for $v$ in \eqref{Th1Eq1} holds true. We now consider the remaining case $q+s=n/\(n-2\)$. In this case, we define
$$\sigma_R:=\frac{n+\ln\(1+R\)^{-1}}{n-2}\,.$$
Similarly as in \eqref{Th1St3Eq14}--\eqref{Th1St3Eq15}, we obtain
\begin{equation}\label{Th1St3Eq19}
\left\|v_R\right\|_{L^{\infty}\(B\(0,4\)\backslash B\(0,2\)\)}\le C\big(\left\|v_R\right\|_{L^{\sigma_R,\infty}\(B\(0,5\)\backslash B\(0,1\)\)}+1\big)
\end{equation}
for some constant $C$ independent of $R$. An easy calculation yields
\begin{multline}
\left\|v_R\right\|_{L^{\sigma_R,\infty}\(B\(0,5\)\backslash B\(0,1\)\)}=\frac{R^{\(\(n-2\)\sigma_R-n\)/\sigma_R}}{\ln\(1+R\)^{1/\(1-s\)}}\left\|v\right\|_{L^{\sigma_R,\infty}\(B\(0,5R\)\backslash B\(0,R\)\)}\\
\le C\ln\(1+R\)^{-1/\(1-s\)}\left\|v\right\|_{L^{\sigma_R,\infty}\(B\(0,5R\)\backslash B\(0,R\)\)}\label{Th1St3Eq20}
\end{multline}
for some constant $C$ independent of $R$. Moreover, it follows from Lemma~\ref{Lem4} that
\begin{align}\label{Th1St3Eq21}
\left\|v\right\|_{L^{\sigma_R,\infty}\(B\(0,5R\)\backslash B\(0,R\)\)}&\le\Lambda_0\big(\(\(n-2\)\sigma_R-n\)^{-1/\(1-s\)}+1\big)\nonumber\\
&=\Lambda_0\big(\ln\(1+R\)^{1/\(1-s\)}+1\big)
\end{align}
By putting together \eqref{Th1St3Eq19}--\eqref{Th1St3Eq21}, we obtain that \eqref{Th1St3Eq17} also holds true in this case. By proceeding as in \eqref{Th1St3Eq18}, we then conclude that the upper bound estimate for $v$ in \eqref{Th1Eq1} holds true in case $s<1$.

\proof[Case $s\ge1$] 
Note that by using \eqref{Eq2}--\eqref{Eq5}, we obtain that $s\ge1$ implies that $q+s>n/\(n-2\)$. In this case, we begin with proving the upper bound estimate for $v$. For any $R\ge1$, we let $v_R$ and $g_R$ be as in \eqref{Th1St3Eq10} and \eqref{Th1St3Eq11}. By applying Step~\ref{Th1St1} and using \eqref{Eq3}, we obtain 
\begin{equation}\label{Th1St3Eq22}
\left|g_R\(y\)\right|\le C\left|y\right|^{-2}\left|v_R\(y\)\right|
\end{equation}
We deduce the upper bound estimate for $v$ from \eqref{Th1St3Eq22} by proceeding as in \eqref{Th1St3Eq4}--\eqref{Th1St3Eq9}. We then define $u_R$ and $f_R$ as in \eqref{Th1St3Eq1} and \eqref{Th1St3Eq2}. By using the upper bound estimate for $v$, together with Step~\ref{Th1St1} and \eqref{Eq5} in case $r\ge1$, we obtain
\begin{align}
\left|f_R\(y\)\right|&\le\left\{\begin{aligned}&CR^{n-\(p+r\)\(n-2\)}\left|y\right|^{-p\(n-2\)}\(\left|u_R\(y\)\right|+1\)&&\text{if }r<1\\
&CR^{-p\(n-2-n/b\)}\left|y\right|^{-2-p\(n-2-n/b\)}\left|u_R\(y\)\right|&&\text{if }r\ge1\end{aligned}\right.\nonumber\\
&\le\left\{\begin{aligned}&C'\left|y\right|^{-p\(n-2\)}\(\left|u_R\(y\)\right|+1\)&&\text{if }r<1\\
&C'\left|y\right|^{-2-p\(n-2-n/b\)}\left|u_R\(y\)\right|&&\text{if }r\ge1\end{aligned}\right.\label{Th1St3Eq23}
\end{align}
for some constants $C,C'$ independent of $y$ and $R$. We then deduce the upper bound estimate for $u$ from \eqref{Th1St3Eq23} by proceeding as in \eqref{Th1St3Eq13}--\eqref{Th1St3Eq18}. This ends the proof of Step~\ref{Th1St3}.
\endproof

The last step in the proof of Theorem~\ref{Th1} is as follows:

\begin{step}\label{Th1St4}
Assume that $u,v\ge0$ in $\R^n$. Then either $u\equiv v\equiv0$ in $\R^n$ or the lower bound estimates \eqref{Th1Eq3} hold true.
\end{step}

\proof[Proof of Step~\ref{Th1St4}]
By applying the strong maximum principle and since $u,v\in L^\infty\(\R^n\)$, we obtain that either $u\equiv v\equiv0$ or $u,v>0$ in $\R^n$. We assume that we are in the latter case. 
By using Step~\ref{Lem4St1}, we obtain
$$u\(x\)\ge\int_{B\(0,1\)}\Gamma\(x-y\)v\(y\)^pu\(y\)^rdy\ge C\(1+\left|x\right|^{n-2}\)^{-1}\qquad\forall x\in\R^n$$
for some constant $C$ independent of $x$, and thus we obtain that the lower bound estimate for $u$ in \eqref{Th1Eq3} holds true. Now we establish the lower bound estimate for $v$. By using Step~\ref{Lem4St1} and the lower bound estimate for $u$, we obtain
\begin{equation}\label{Th1St4Eq10} 
\int_{\R^n}\left|x-y\right|^{2-n}\(1+\left|y\right|^{n-2}\)^{-q}v\(y\)^sdy\le Cv\(x\)\qquad\forall x\in\R^n
\end{equation}
for some constant $C$ independent of $x$. We define
$$A_0:=\inf_{\left|x\right|<1}v\(x\),\quad A_k:=\inf_{2^{k-1}<\left|x\right|<2^k}v\(x\)\quad\forall k\in\N\backslash\left\{0\right\},$$
$$I_{0,k}:=\inf_{2^{k-1}<\left|x\right|<2^k}\int_{B\(0,1\)}\left|x-y\right|^{2-n}dy\quad\forall k\in\N\backslash\left\{0\right\},$$
and
$$I_{j,k}:=\inf_{2^{k-1}<\left|x\right|<2^k}\int_{B\(0,2^j\)\backslash B\(0,2^{j-1}\)}\left|x-y\right|^{2-n}dy\quad\forall j,k\in\N\backslash\left\{0\right\}.$$
It follows from \eqref{Th1St4Eq10} that
\begin{equation}\label{Th1St4Eq11} 
\sum_{j=0}^\infty2^{-jq\(n-2\)}A_j^sI_{j,k}\le CA_k.
\end{equation}
for some constant $C$ independent of $j$ and $k$. An easy calculation gives
\begin{equation}\label{Th1St4Eq12} 
I_{j,k}\ge c2^{nj-\(n-2\)k}\qquad\forall j\in\left\{0,\dotsc,k\right\}
\end{equation}
for some constant $c>0$ independent of $k$. It follows from \eqref{Th1St4Eq11} and \eqref{Th1St4Eq12} that
\begin{equation}\label{Th1St4Eq13}
A_k\ge\left\{\begin{aligned}&c2^{-k\(n-2\)}&&\text{if }q+s>n/\(n-2\)\\&c2^{-k\(n-2\)}k^{1/\(1-s\)}&&\text{if }q+s=n/\(n-2\)\\&c2^{-k\(\(n-2\)q-2\)/\(1-s\)}&&\text{if }q+s<n/\(n-2\)\end{aligned}\right.
\end{equation}
for some constant $c>0$ independent of $k$. It follows from \eqref{Th1St4Eq13} that the lower bound estimate for $v$ in \eqref{Th1Eq3} holds true. 
This ends the proof of Step~\ref{Th1St4} and also Theorem~\ref{Th1}.
\endproof

\section{Symmetry results}\label{Sec4}

This section is devoted to the proofs of Theorems~\ref{Th2}--\ref{Th4}. We assume that \eqref{Eq2} and \eqref{Eq5} hold true. We fix a solution $\(u,v\)\in L^a\(\R^n\)\times L^b\(\R^n\)$ of \eqref{Eq1}. By applying the strong maximum principle and since $u,v\in L^\infty\(\R^n\)$, we obtain that either $u\equiv v\equiv0$ or $u,v>0$ in $\R^n$. We assume that we are in the latter case. In case $q+s<n/\(n-2\)$, we assume moreover that \eqref{Th3Eq} holds true.

\smallskip
This proof relies on the moving plane method and follows in great part the lines of Liu and Ma~\cite{LiuMa1} but with the key difference in case $q+s<n/\(n-2\)$ that we consider different powers in the definitions of $\overline{U}_\lambda$ and $\overline{V}_\lambda$ below. For any $\lambda>0$, we define
$$\Sigma_\lambda:=\left\{x=\(x_1,x'\)\in\R^n:\,x_1<\lambda\right\}$$
and
$$U_\lambda\(x\):=u_\lambda\(x\)-u\(x\)\quad\text{and}\quad V_\lambda\(x\):=v_\lambda\(x\)-v\(x\)$$
where
$$u_\lambda\(x\):=u\(x_\lambda\),\quad v_\lambda\(x\):=v\(x_\lambda\),\quad\text{and}\quad x_\lambda:=\(2\lambda-x_1,x'\).$$
By using \eqref{Eq1} and applying the mean value theorem, we obtain
\begin{equation}\label{Th2Eq1}
-\Delta U_\lambda=rv^p\xi_{u,u}^{r-1}U_\lambda+p\xi_{u,v}^{p-1}u_\lambda^rV_\lambda\qquad\text{in }\R^n
\end{equation}
and
\begin{equation}\label{Th2Eq2}
-\Delta V_\lambda=su^q\xi_{v,v}^{s-1}V_\lambda+q\xi_{v,u}^{q-1}v_\lambda^sU_\lambda\qquad\text{in }\R^n
\end{equation}
for some functions $\xi_{u,u},\xi_{u,v},\xi_{v,u},\xi_{v,v}:\R^n\to\R$ such that
\begin{equation}\label{Th2Eq3} 
\min\(u,u_\lambda\)\le\xi_{u,u},\xi_{v,u}\le\max\(u,u_\lambda\)
\end{equation}
and
\begin{equation}\label{Th2Eq4} 
\min\(v,v_\lambda\)\le\xi_{u,v},\xi_{v,v}\le\max\(v,v_\lambda\).
\end{equation}
We define
$$\left\{\begin{aligned}&\theta_u=\theta_v:=\frac{n-2}{2}&&\text{if }q+s>n/\(n-2\)\\&\theta_u=\theta_v+\varepsilon:=\frac{n-2}{2}&&\text{if }q+s=n/\(n-2\)\\&\left\{\begin{aligned}&\theta_u:=\theta_v+\frac{n-\(n-2\)\(q+s\)}{1-s}+\varepsilon\\&\theta_v:=\min\Big(\frac{n-2}{2},\frac{\(n-2\)q-2}{1-s}-2\varepsilon\Big)\end{aligned}\right.&&\text{if }q+s<n/\(n-2\)\end{aligned}\right.$$
where $\varepsilon$ is a small positive parameter to be fixed later on. One can easily see that if $\varepsilon$ is chosen small enough, then
\begin{equation}\label{Th2Eq5} 
0<\theta_u<n-2
\end{equation}
and
\begin{equation}\label{Th2Eq6}
0<\theta_v<\left\{\begin{aligned}&n-2&&\text{if }q+s\ge n/\(n-2\)\\&\frac{\(n-2\)q-2}{1-s}&&\text{if }q+s<n/\(n-2\).\end{aligned}\right.
\end{equation}
We then define
$$\overline{U}_\lambda\(x\):=\(1+\left|x\right|^2\)^{\theta_u/2}U_\lambda\(x\)\quad\text{and}\quad\overline{V}_\lambda\(x\):=\(1+\left|x\right|^2\)^{\theta_v/2}V_\lambda\(x\).$$
It follows from \eqref{Th2Eq5}, \eqref{Th2Eq6}, and Theorem~\ref{Th1} that
\begin{equation}\label{Th2Eq7} 
\lim_{\left|x\right|\to\infty}\overline{U}_\lambda\(x\)=0\quad\text{and}\quad\lim_{\left|x\right|\to\infty}\overline{V}_\lambda\(x\)=0.
\end{equation}
Moreover, it follows from \eqref{Th2Eq1} and \eqref{Th2Eq2} that
\begin{multline}\label{Th2Eq8}
-\Delta\overline{U}_\lambda+\frac{2\theta_u}{1+\left|x\right|^2}\<x,\nabla\overline{U}_\lambda\>+\frac{\theta_u\(n+\(n-2-\theta_u\)\left|x\right|^2\)}{\(1+\left|x\right|^2\)^2}\overline{U}_\lambda\\
=rv^p\xi_{u,u}^{r-1}\overline{U}_\lambda+p\xi_{u,v}^{p-1}u_\lambda^r\(1+\left|x\right|^2\)^{\(\theta_u-\theta_v\)/2}\overline{V}_\lambda\qquad\text{in }\R^n
\end{multline}
and
\begin{multline}\label{Th2Eq9}
-\Delta\overline{V}_\lambda+\frac{2\theta_v}{1+\left|x\right|^2}\<x,\nabla\overline{V}_\lambda\>+\frac{\theta_v\(n+\(n-2-\theta_v\)\left|x\right|^2\)}{\(1+\left|x\right|^2\)^2}\overline{V}_\lambda\\
=su^q\xi_{v,v}^{s-1}\overline{V}_\lambda+q\xi_{v,u}^{q-1}v_\lambda^s\(1+\left|x\right|^2\)^{\(\theta_v-\theta_u\)/2}\overline{U}_\lambda\qquad\text{in }\R^n.
\end{multline}

\smallskip
As a first step, we prove the following result:

\begin{step}\label{Th2St1}
Let $\varepsilon$ be the positive constant in the definition of $\theta_u$ and $\theta_v$. There exists $\delta_0,\varepsilon_0>0$ such that if $\varepsilon\in\(0,\varepsilon_0\)$, then for any $\delta\in\(0,\delta_0\)$ and $\lambda\in\R$, there exists $R_{\delta,\varepsilon,\lambda}>0$ such that
\begin{align}
&pv\(x\)^{p-1}u\(x\)^r<\frac{\theta_u\(n+\(n-2-\theta_u\)\left|x\right|^2\)}{2\(1+\left|x\right|^2\)^{2+\(\theta_u-\theta_v\)/2}}\,,\label{Th2St1Eq1}\allowdisplaybreaks\\
&rv\(x\)^pu\(x\)^{r-1}<\frac{\theta_u\(n+\(n-2-\theta_u\)\left|x\right|^2\)}{2\(1+\left|x\right|^2\)^2}&\text{if }r\ge1,\label{Th2St1Eq2}\allowdisplaybreaks\\
&rv\(x\)^pu_\mu\(x\)^{r-1}<\frac{\theta_u\(n+\(n-2-\theta_u\)\left|x\right|^2\)}{2\(1+\left|x\right|^2\)^2}&\text{if }r<1,\label{Th2St1Eq3}\allowdisplaybreaks\\
&qu\(x\)^{q-1}v\(x\)^s<\frac{\delta\theta_v\(n+\(n-2-\theta_v\)\left|x\right|^2\)}{\(1+\left|x\right|^2\)^{2+\(\theta_v-\theta_u\)/2}}\,,\label{Th2St1Eq4}\allowdisplaybreaks\\
&su\(x\)^qv\(x\)^{s-1}<\frac{\(1-\delta\)\theta_v\(n+\(n-2-\theta_v\)\left|x\right|^2\)}{\(1+\left|x\right|^2\)^2}&\text{if }s\ge1,\label{Th2St1Eq5}\allowdisplaybreaks\\
&su\(x\)^qv_\mu\(x\)^{s-1}<\frac{\(1-\delta\)\theta_v\(n+\(n-2-\theta_v\)\left|x\right|^2\)}{\(1+\left|x\right|^2\)^2}&\text{if }s<1\label{Th2St1Eq6}
\end{align}
for all $\mu\le\lambda$ and $x\in\Sigma_\mu\backslash B\(0,R_{\delta,\varepsilon,\lambda}\)$. 
\end{step}

\proof[Proof of Step~\ref{Th2St1}]
By using \eqref{Eq2}--\eqref{Eq5} and Theorem~\ref{Th1}, we obtain
\begin{align}
&v\(x\)^{p-1}u\(x\)^r\nonumber\\
&\quad=\left\{\begin{aligned}&\OO\big(\left|x\right|^{-\(n-2\)\(p+r-1\)}\big)&&\text{if }q+s>n/\(n-2\)\\&\OO\big(\left|x\right|^{-\(n-2\)\(p+r-1\)}\ln\(1+\left|x\right|\)^{\(p-1\)/\(1-s\)}\big)&&\text{if }q+s=n/\(n-2\)\\&\OO\big(\left|x\right|^{-\(\(n-2\)\(q\(p-1\)+r\(1-s\)\)-2\(p-1\)\)/\(1-s\)}\big)&&\text{if }q+s<n/\(n-2\)\end{aligned}\right.\nonumber\\
&\quad=\oo\big(\left|x\right|^{-2+\theta_v-\theta_u}\big)\label{Th2St1Eq7}
\end{align}
as $\left|x\right|\to\infty$ provided $\varepsilon$ is small enough so that
$$\varepsilon<\left\{\begin{aligned}&\(n-2\)\(p+r\)-n&&\text{if }q+s=n/\(n-2\)\\&2\(\(n-2\)a-n\)\(p-s+1\)/\(n\(1-s\)\)&&\text{if }q+s<n/\(n-2\).\end{aligned}\right.$$
Moreover, by using \eqref{Eq2}--\eqref{Eq5} and Theorem~\ref{Th1}, we obtain
\begin{align}
&v\(x\)^pu\(x\)^{r-1}\nonumber\\
&\quad=\left\{\begin{aligned}&\OO\big(\left|x\right|^{-\(n-2\)\(p+r-1\)}\big)&&\text{if }q+s>n/\(n-2\)\\&\OO\big(\left|x\right|^{-\(n-2\)\(p+r-1\)}\ln\(1+\left|x\right|\)^{\(r-1\)/\(1-s\)}\big)&&\text{if }q+s=n/\(n-2\)\\&\OO\big(\left|x\right|^{-\(\(n-2\)\(pq-\(r-1\)\(s-1\)\)-2p\)/\(1-s\)}\big)&&\text{if }q+s<n/\(n-2\)\end{aligned}\right.\nonumber\\
&\quad=\oo\big(\left|x\right|^{-2}\big),\label{Th2St1Eq8}
\end{align}
\begin{align}
u\(x\)^{q-1}v\(x\)^s&=\left\{\begin{aligned}&\OO\big(\left|x\right|^{-\(n-2\)\(q+s-1\)}\big)&&\text{if }q+s>n/\(n-2\)\\&\OO\big(\left|x\right|^{-2}\ln\(1+\left|x\right|\)^{s/\(1-s\)}\big)&&\text{if }q+s=n/\(n-2\)\\&\OO\big(\left|x\right|^{-\(\(n-2\)\(q+s-1\)-2s\)/\(1-s\)}\big)&&\text{if }q+s<n/\(n-2\)\end{aligned}\right.\nonumber\\
&=\oo\big(\left|x\right|^{-2+\theta_u-\theta_v}\big),\label{Th2St1Eq9}
\end{align}
and
\begin{align}
u\(x\)^qv\(x\)^{s-1}&=\left\{\begin{aligned}&\OO\big(\left|x\right|^{-\(n-2\)\(q+s-1\)}\big)&&\text{if }q+s>n/\(n-2\)\\&\OO\big(\left|x\right|^{-2}\ln\(1+\left|x\right|\)^{-1}\big)&&\text{if }q+s=n/\(n-2\)\\&\OO\big(\left|x\right|^{-2}§\big)&&\text{if }q+s<n/\(n-2\)\end{aligned}\right.\nonumber\\
&=\oo\big(\left|x\right|^{-2}\big)\quad\text{if }q+s\ge n/\(n-2\)\label{Th2St1Eq10}
\end{align}
as $\left|x\right|\to\infty$. Note that by using \eqref{Eq2}--\eqref{Eq5}, we obtain that $s\ge1$ implies that $q+s>n/\(n-2\)$. The estimates \eqref{Th2St1Eq1}, \eqref{Th2St1Eq2}, \eqref{Th2St1Eq4}, and \eqref{Th2St1Eq5} then follow from \eqref{Th2St1Eq7}--\eqref{Th2St1Eq10}. Now we prove the remaining estimates \eqref{Th2St1Eq3} and \eqref{Th2St1Eq5}. By applying Young's inequality, we obtain that for any $\delta>0$, there exists $C_\delta>0$ such that
\begin{equation}\label{Th2St1Eq11} 
\left|x_\mu\right|^2=\left|x\right|^2+4\mu\(\mu-x_1\)\le\(1+\delta\)\left|x\right|^2+C_\delta\lambda^2.
\end{equation}
for all $\mu\le\lambda$ and $x\in\Sigma_\mu$. The estimate \eqref{Th2St1Eq3} then follows from \eqref{Th2St1Eq8}, \eqref{Th2St1Eq11}, and Theorem~\ref{Th1}. In case $q+s\ge n/\(n-2\)$, the estimate \eqref{Th2St1Eq5} follows from \eqref{Th2St1Eq10}, \eqref{Th2St1Eq11}, and Theorem~\ref{Th1}. In case $q+s<n/\(n-2\)$, the estimate \eqref{Th2St1Eq5} follows from \eqref{Th3Eq} and \eqref{Th2St1Eq11} provided $\delta$ and $\varepsilon$ are small enough. This ends the proof of Step~\ref{Th2St1}. 
\endproof

Remark that in order to have \eqref{Th2St1Eq6} in case $q+s<n/\(n-2\)$ and $s\ne0$, it is necessary to make a decay assumption such as \eqref{Th3Eq}. Moreover, the best constant in \eqref{Th3Eq} with this approach is given by maximizing $\theta_v\(n-2-\theta_v\)/s$ with $\varepsilon=0$ under the constraint \eqref{Th2Eq6} which gives the constant $C_{n,q,s}$ defined in Theorem~\ref{Th3}.

\smallskip
Then we prove the following result:

\begin{step}\label{Th2St2}
Assume that $\varepsilon\in\(0,\varepsilon_0\)$ and let $\delta\in\(0,\delta_0\)$ where $\delta_0$ and $\varepsilon_0$ are as in Step~\ref{Th2St1}. Then $U_\lambda\ge0$ and $V_\lambda\ge0$ in $\Sigma_\lambda$ for all $\lambda\le-R_{\delta,\varepsilon,0}$ where $R_{\delta,\varepsilon,0}$ is as in Step~\ref{Th2St1}.
\end{step}

\proof[Proof of Step~\ref{Th2St2}]
We assume by contradiction that there exists $\lambda\le-R_{\delta,\varepsilon,0}$ and $y_0\in\Sigma_\lambda$ such that either $U_\lambda\(y_0\)<0$ or $V_\lambda\(y_0\)<0$. Since the proof is similar in cases $U_\lambda\(y_0\)<0$ and $V_\lambda\(y_0\)<0$, we will assume that $U_\lambda\(y_0\)<0$.  In this case, by using \eqref{Th2Eq7} and since $\overline{U}_\lambda\equiv0$ on $\partial\Sigma_\lambda$, we obtain that there exists $y_1\in\Sigma_\lambda$ such that 
\begin{equation}\label{Th2St2Eq1}
\overline{U}_\lambda\(y_1\):=\min\left\{\overline{U}_\lambda\(y\):\,y\in\Sigma_\lambda\right\}<0.
\end{equation}
It follows from \eqref{Th2Eq8} and \eqref{Th2St2Eq1} that
\begin{multline}\label{Th2St2Eq2}
p\xi_{u,v}\(y_1\)^{p-1}u_\lambda\(y_1\)^r\(1+\left|y_1\right|^2\)^{\(\theta_u-\theta_v\)/2}\overline{V}_\lambda\(y_1\)\\
\le\bigg(\frac{\theta_u\(n+\(n-2-\theta_u\)\left|y_1\right|^2\)}{\(1+\left|y_1\right|^2\)^2}-rv\(y_1\)^p\xi_{u,u}\(y_1\)^{r-1}\bigg)\overline{U}\(y_1\).
\end{multline}
Since $\left|y_1\right|=-y_1>-\lambda\ge R_{\delta,\varepsilon,0}$ and $U_\lambda\(y_1\)<0$, it follows from \eqref{Th2Eq3}, \eqref{Th2St1Eq2}, and \eqref{Th2St1Eq3} that
\begin{equation}\label{Th2St2Eq}
rv\(y_1\)^p\xi_{u,u}\(y_1\)^{r-1}<\frac{\theta_u\(n+\(n-2-\theta_u\)\left|y_1\right|^2\)}{2\(1+\left|y_1\right|^2\)^2}\,.
\end{equation}
By using \eqref{Th2Eq4}, \eqref{Th2St2Eq1}, \eqref{Th2St2Eq2}, and \eqref{Th2St2Eq} and since $p\ge1$ and $r\ge0$, we obtain
\begin{align}\label{Th2St2Eq3}
&pv\(y_1\)^{p-1}u\(y_1\)^r\(1+\left|y_1\right|^2\)^{\(\theta_u-\theta_v\)/2}\overline{V}_\lambda\(y_1\)\nonumber\\
&\qquad<\frac{\theta_u\(n+\(n-2-\theta_u\)\left|y_1\right|^2\)}{2\(1+\left|y_1\right|^2\)^2}\overline{U}_\lambda\(y_1\)<0.
\end{align}
By using \eqref{Th2Eq7} and since $\overline{V}_\lambda\equiv0$ on $\partial\Sigma_\lambda$, it follows that there exists $y_2\in\Sigma_\lambda$ such that 
\begin{equation}\label{Th2St2Eq4}
\overline{V}_\lambda\(y_2\):=\min\left\{\overline{V}_\lambda\(y\):\,y\in\Sigma_\lambda\right\}<0.
\end{equation}
By repeating the above arguments and using \eqref{Th2Eq3}, \eqref{Th2Eq4}, \eqref{Th2Eq9}, \eqref{Th2St1Eq5}, and \eqref{Th2St1Eq6}, we obtain
\begin{align}\label{Th2St2Eq5}
&qu\(y_2\)^{q-1}v\(y_2\)^s\(1+\left|y_2\right|^2\)^{\(\theta_v-\theta_u\)/2}\overline{U}_\lambda\(y_2\)\nonumber\\
&\qquad<\frac{\delta\theta_v\(n+\(n-2-\theta_v\)\left|y_2\right|^2\)}{\(1+\left|y_2\right|^2\)^2}\overline{V}_\lambda\(y_2\)<0.
\end{align}
It follows from \eqref{Th2St2Eq1}, \eqref{Th2St2Eq4}, and \eqref{Th2St2Eq5} that
\begin{align}\label{Th2St2Eq6}
&qu\(y_2\)^{q-1}v\(y_2\)^s\(1+\left|y_2\right|^2\)^{\(\theta_v-\theta_u\)/2}\overline{U}_\lambda\(y_1\)\nonumber\\
&\qquad<\frac{\delta\theta_v\(n+\(n-2-\theta_v\)\left|y_2\right|^2\)}{\(1+\left|y_2\right|^2\)^2}\overline{V}_\lambda\(y_1\)<0.
\end{align}
It follows from \eqref{Th2St2Eq3} and \eqref{Th2St2Eq6} that
\begin{multline}\label{Th2St2Eq7}
\big(pv\(y_1\)^{p-1}u\(y_1\)^r\big)\cdot\big(qu\(y_2\)^{q-1}v\(y_2\)^s\big)\\
>\frac{\theta_u\(n+\(n-2-\theta_u\)\left|y_1\right|^2\)}{2\(1+\left|y_1\right|^2\)^{2+\(\theta_u-\theta_v\)/2}}\cdot\frac{\delta\theta_v\(n+\(n-2-\theta_v\)\left|y_2\right|^2\)}{\(1+\left|y_2\right|^2\)^{2+\(\theta_v-\theta_u\)/2}}\,.
\end{multline}
We obtain a contradiction by putting together \eqref{Th2St1Eq1}, \eqref{Th2St1Eq4}, and \eqref{Th2St2Eq7}. This ends the proof of Step~\ref{Th2St1}.
\endproof

We define
\begin{equation}\label{Th2Eq10}
\lambda_0:=\sup\left\{\lambda\in\R:\,\forall\mu\le\lambda,\,U_\mu\ge0\text{ and }V_\mu\ge0\text{ in }\Sigma_\mu\right\}.
\end{equation}
It follows from Step~\ref{Th2St2} that if $\varepsilon\in\(0,\varepsilon_0\)$, then $\lambda_0>-\infty$. On the other hand, it follows from Theorem~\ref{Th1} that $U_\lambda\(0\)<0$ for large $\lambda>0$ hence we obtain $\lambda_0<+\infty$.

\smallskip
The last step in the proof of Theorems~\ref{Th2} and~\ref{Th3} is as follows:

\begin{step}\label{Th2St3}
Assume that $\varepsilon\in\(0,\varepsilon_0\)$ where $\varepsilon_0$ is as in Step~\ref{Th2St1}. Let $\lambda_0$ be as in \eqref{Th2Eq10}. Then $U_{\lambda_0}\equiv0$ and $V_{\lambda_0}\equiv0$ in $\Sigma_{\lambda_0}$.
\end{step}

\proof[Proof of Step~\ref{Th2St3}]
It follows from the definition of $\lambda_0$ and the continuity of $U_\mu$ and $V_\mu$ with respect to $\mu$ that $U_{\lambda_0}\ge0$ and $V_{\lambda_0}\ge0$ in $\Sigma_{\lambda_0}$. By using \eqref{Th2Eq1} and \eqref{Th2Eq2} and applying the strong maximum principle, we then obtain that either $U_{\lambda_0}\equiv V_{\lambda_0}\equiv0$ or $U_{\lambda_0},V_{\lambda_0}>0$ in $\Sigma_{\lambda_0}$. We assume by contradiction that $U_{\lambda_0},V_{\lambda_0}>0$ in $\Sigma_{\lambda_0}$. By definition of $\lambda_0$, we obtain that there exists a sequence of real numbers $\(\mu_\alpha\)_{\alpha\in\N}$ such that $\mu_\alpha\searrow\lambda_0$ as $\alpha\to\infty$ and either
\begin{equation}\label{Th2St3Eq1}
\inf\left\{U_{\mu_\alpha}\(y\):\,y\in\Sigma_{\mu_\alpha}\right\}<0\quad\text{or}\quad\inf\left\{V_{\mu_\alpha}\(y\):\,y\in\Sigma_{\mu_\alpha}\right\}<0
\end{equation}
for all $\alpha$. Since the proof is similar in both cases, we will assume that the first inequality in \eqref{Th2St3Eq1} holds true. By using \eqref{Th2Eq7} and since $\overline{U}_\lambda\equiv0$ on $\partial\Sigma_\lambda$, we then obtain that there exists $y_\alpha\in\Sigma_{\mu_\alpha}$ such that 
\begin{equation}\label{Th2St3Eq2}
\overline{U}_{\mu_\alpha}\(y_\alpha\)=\min\left\{\overline{U}_{\mu_\alpha}\(y\):\,y\in\Sigma_{\mu_\alpha}\right\}<0.
\end{equation}
By extracting a subsequence, we can assume that either $\left|y_\alpha\right|>R_{\delta,\varepsilon,\mu_0}$ or $\left|y_\alpha\right|\le R_{\delta,\varepsilon,\mu_0}$ for all $\alpha$ where $\delta\in\(0,\delta_0\)$ is fixed and $\delta_0$ and $R_{\delta,\varepsilon,\mu_0}$ are as in Step~\ref{Th2St1}. In case $\left|y_\alpha\right|\le R_{\delta,\varepsilon,\mu_0}$, since $y_\alpha\in\Sigma_{\mu_\alpha}$ and $\mu_\alpha\to\lambda_0$ as $\alpha\to\infty$, we obtain that up to a subsequence $y_\alpha\to y_\infty\in\overline{\Sigma_{\lambda_0}}$ as $\alpha\to\infty$. By passing to the limit into \eqref{Th2St3Eq2} and since $U_{\lambda_0}>0$ in $\Sigma_{\lambda_0}$, we obtain
\begin{equation}\label{Th2St3Eq3}
\overline{U}_{\lambda_0}\(y_\infty\)=\min\left\{\overline{U}_{\lambda_0}\(y\):\,y\in\overline{\Sigma_{\lambda_0}}\right\}=0\,\text{ and }\,\nabla\overline{U}_{\lambda_0}\(y_\infty\)=0.
\end{equation}
It follows from \eqref{Th2St3Eq3} and the definition of $\overline{U}_{\lambda_0}$ that
\begin{equation}\label{Th2St3Eq4}
U_{\lambda_0}\(y_\infty\)=\min\left\{U_{\lambda_0}\(y\):\,y\in\overline{\Sigma_{\lambda_0}}\right\}=0\,\text{ and }\,\nabla U_{\lambda_0}\(y_\infty\)=0.
\end{equation}
Since $U_{\lambda_0}$ is a solution of \eqref{Th2Eq1} and $U_{\lambda_0},V_{\lambda_0}>0$ in $\Sigma_{\lambda_0}$, we obtain a contradiction between \eqref{Th2St3Eq3} and Hopf's Lemma. Now we consider the case $\left|y_\alpha\right|>R_{\delta,\varepsilon,\mu_0}$. Since $\mu_\alpha<\mu_0$, by proceeding as in \eqref{Th2St2Eq2}--\eqref{Th2St2Eq4}, we obtain that there exists $z_\alpha\in\Sigma_{\mu_\alpha}$ such that 
$$\overline{V}_{\mu_\alpha}\(z_\alpha\)=\min\left\{\overline{V}_{\mu_\alpha}\(y\):\,y\in\Sigma_{\mu_\alpha}\right\}<0.$$
In case $\left|z_\alpha\right|>R_{\delta,\varepsilon,\mu_0}$, we obtain a contradiction by repeating the arguments in \eqref{Th2St2Eq5}--\eqref{Th2St2Eq7}. In case $\left|z_\alpha\right|\le R_{\delta,\varepsilon,\mu_0}$, we also obtain a contradiction by proceeding as in case $\left|y_\alpha\right|\le R_{\delta,\varepsilon,\mu_0}$. This ends the proof of Step~\ref{Th2St3}.
\endproof

\proof[End of proof of Theorems~\ref{Th2} and~\ref{Th3}] 
Since we can apply the moving plane method in any direction, it follows from Step~\ref{Th2St3} that $u$ and $v$ are radially symmetric and decreasing about some point $x_0\in\R^n$. This ends the proof of Theorem~\ref{Th2} and~\ref{Th3}.
\endproof

Finally we prove Theorem~\ref{Th4}.

\proof[Proof of Theorem~\ref{Th4}]
We assume that $q+s<n/\(n-2\)$ and \eqref{Eq10} and \eqref{Eq11} hold true. By applying Step~\ref{Lem4St1}, we obtain
\begin{equation}\label{Th4Eq1}
v\(x\)=\Gamma\ast\(u^qv^s\)\(x\)=I_1\(x\)+I_2\(x\)
\end{equation}
where
$$I_1\(x\):=\int_{B\(0,\sqrt{\left|x\right|}\)}\Gamma\(x-y\)u\(y\)^qv\(y\)^sdy$$
and
$$I_2\(x\):=\int_{\R^n\backslash B\(0,\sqrt{\left|x\right|}\)}\Gamma\(x-y\)u\(y\)^qv\(y\)^sdy.$$
By using \eqref{Eq10} and \eqref{Eq11}, we obtain
\begin{align}
I_1\(x\)&=\oo\(\int_{B\(0,\left|x\right|\)}\left|x-y\right|^{2-n}\left|y\right|^{\(2s-q\(n-2\)\)/\(1-s\)}dy\)\nonumber\\
&=\oo\big(\left|x\right|^{\(2-q\(n-2\)\)/\(1-s\)}\big)=\oo\(v\(x\)\)\label{Th4Eq2}
\end{align}
as $\left|x\right|\to\infty$. It follows from \eqref{Eq10}, \eqref{Eq11},\eqref{Th4Eq1}, and \eqref{Th4Eq2} that
\begin{align}
\frac{1}{\ell_u^q\ell_v^{s-1}}&=\lim_{\left|x\right|\to\infty}\(\int_{\R^n\backslash B\(0,\sqrt{\left|x\right|}\)}\Gamma\(x-y\)\frac{\left|x\right|^{\(q\(n-2\)-2\)/\(1-s\)}}{\left|y\right|^{\(q\(n-2\)-2s\)/\(1-s\)}}dy\)\nonumber\\
&=\lim_{\left|x\right|\to\infty}\(\int_{\R^n\backslash B\(0,1/\sqrt{\left|x\right|}\)}\Gamma\Big(\frac{x}{\left|x\right|}-y\Big)\left|y\right|^{\(2s-q\(n-2\)\)/\(1-s\)}dy\)\nonumber\\
&=\int_{\R^n}\Gamma\(\omega-y\)\left|y\right|^{\(2s-q\(n-2\)\)/\(1-s\)}dy\label{Th4Eq3}
\end{align}
where $\omega$ is an arbitrary point on $\partial B\(0,1\)$. By observing that the function $\varphi:\R^n\backslash\left\{0\right\}\to\R$ defined as 
$$\varphi\(x\)=\frac{\(1-s\)^2\left|x\right|^{\(2-q\(n-2\)\)/\(1-s\)}}{\(\(n-2\)q-2\)\(n-\(n-2\)\(q+s\)\)}$$
is a distributional solution of the equation
$$-\Delta\varphi=\left|x\right|^{\(2s-q\(n-2\)\)/\(1-s\)}\qquad\text{in }\R^n$$
and $\varphi\(x\)\to0$ as $\left|x\right|\to\infty$, we obtain
\begin{multline}
\int_{\R^n}\Gamma\(\omega-y\)\left|y\right|^{\(2s-q\(n-2\)\)/\(1-s\)}dy\\
=\varphi\(\omega\)=\frac{\(1-s\)^2}{\(\(n-2\)q-2\)\(n-\(n-2\)\(q+s\)\)}\,.\label{Th4Eq4}
\end{multline}
It follows from \eqref{Eq10}, \eqref{Eq11},\eqref{Th4Eq3}, and \eqref{Th4Eq4} that
$$c_u^qc_v^{s-1}=\frac{\(\(n-2\)q-2\)\(n-\(n-2\)\(q+s\)\)}{\(1-s\)^2}\,.$$
This ends the proof of Theorem~\ref{Th4}.
\endproof

\end{document}